\documentclass[12pt]{article}

\usepackage{graphicx}
\usepackage{enumerate,csquotes,ulem}

\makeatletter
\newcommand{\labitem}[2]{%
\def\@itemlabel{\textbf{#1}}
\item
\def\@currentlabel{#1}\label{#2}}
\makeatother

\usepackage{natbib}
\usepackage{url} 

\usepackage{graphicx,bm,url,microtype}
\usepackage{hyperref}
\usepackage{paralist}
\usepackage{diagbox}

\usepackage{authblk}
\usepackage{amsfonts,amssymb,amsthm,amsmath, amsthm}
\usepackage{natbib}
\usepackage[colorinlistoftodos]{todonotes}
\usepackage[utf8]{inputenc}
\newtheorem{proposition}{{\sc\bf Proposition}}
\newtheorem{theorem}{{\sc\bf Theorem}}
\newtheorem{lemma}{{\sc\bf Lemma}}

\newtheorem{example}{Example}

\newtheorem{corollary}{{\sc\bf Corollary}}

\definecolor{violet}{rgb}{0.7,0,0.6}

\newcommand{\cd}{\ensuremath{\rightsquigarrow}}

\newcommand{\R}{\ensuremath{\mathbb{R}}}
\newcommand{\F}{\ensuremath{\mathcal{F}}}

\def \E{\text{\rm E}}
\def \P{\text{\rm P}}

\DeclareMathOperator*{\argmin}{arg\,min}

\begin{document}
	\begin{center}
		\Large \bf On uniqueness of the set of $k$-means
	\end{center}
	\footnotesize
	\begin{center}
		Javier Cárcamo$^1$, Antonio Cuevas$^2$ and Luis A. Rodríguez$^3$\\
		$^1$ Departamento de Matem\'aticas, Universidad del País Vasco/Euskal Herriko Unibertsitatea\\
		$^2$ Departamento de Matem\'aticas, Universidad Aut\'onoma de Madrid\\ and Instituto de Ciencias Matem\'aticas ICMAT (CSIC-UAM-UCM-UC3M)\\
		$^3$ Institut f\"{u}r Mathematische Stochastik, Georg-August-Universität G\"{o}ttingen
	\end{center}

\normalsize
\begin{abstract}
We provide necessary and sufficient conditions for the uniqueness of the $k$-means set of a probability distribution. This uniqueness problem is related to the choice of $k$: depending on the underlying distribution, some values of this parameter could lead to multiple sets of $k$-means, which hampers the interpretation of the results and/or the stability of the algorithms.
We give a general assessment on consistency of the empirical $k$-means adapted to the setting of non-uniqueness and determine the asymptotic distribution of the within cluster sum of squares (WCSS).
We also provide statistical characterizations of $k$-means uniqueness in terms of the asymptotic behavior of the empirical WCSS. As a consequence, we derive a bootstrap test for uniqueness of the set of $k$-means. The results are illustrated with examples of different types of non-uniqueness and we check by simulations the performance of the proposed methodology.
\end{abstract}

 \section{Introduction}

The $k$-means procedure is one of the most commonly used techniques for finding a given number of groups in a data set. The notion of $k$-means is a natural idea with a clear interpretation and a great number of relevant applications. However, despite its simplicity, the underlying methodology still has some extraordinarily complex challenges associated with it, such as the choice of the parameter $k$, and several intriguing theoretical and computational aspects. In this paper we focus on the problem of uniqueness of the set of $k$-means, both from a theoretical and practical point of view. As we discuss below, determining whether the set of $k$-means is unique or not could potentially shed light on possible reasonable choices of the parameter $k$, or, at least, to avoid bad selections of $k$.

Given a random element $X$ taking values in a separable Banach space ${\mathcal B}$ and a prefixed natural number $k$, the goal of $k$-means is to find the optimal set of centers of $k$ groups, say $\mu_1,\dots,\mu_k$, such that the expected square distance from $X$ to its nearest center is minimal. More formally, if $X$ induces the probability measure $\P$ on $\mathcal{B}$ (with norm $\Vert \cdot \Vert_{\mathcal{B}}$) and $k\in{\mathbb N}$, the \textit{principal points} or  \textit{$k$-means set} of $\P$ is any subset (of cardinality $k$) of $\mathcal{B}$
minimizing, over all possible sets $\{a_1,\ldots,a_k\}\subset\mathcal{B}$, the quantity
  \begin{equation} \label{kmeans-pop}
  	\Phi(a_1,\ldots,a_k) =  \Phi(\P;a_1,\ldots,a_k)=\E_{\P}\left(
  	\min_{i=1,...,k} \|X - a_i\|_{\mathcal{B}}^2 \right),
  \end{equation}
where `$\E_\P$' stands for the mathematical expectation with respect to $\P$. The (optimal) population \textit{within cluster sum of squares} (WCSS in short) is defined by
\begin{equation}\label{WCSS-population}
 W(k) = W(\P;k)= \inf_{\{a_1,\dots,a_k\}\subset\mathcal{B}} \Phi(a_1,\ldots,a_k).
\end{equation}
Thus, the value $W(k)$ is achieved for any set of $k$-means.

The associated empirical version, based on $n$ independent observations $X_1,\ldots,X_n$ drawn from $X$, corresponds to minimize 
the expression \eqref{kmeans-pop} for the \textit{empirical measure}
$\mathbb{P}_{n}=\frac{1}{n}\,\sum_{i=1}^{n}\delta_{X_{i}}$,
where $\delta_{a}$ stands for the unit point mass at $a$. In other words, we have to minimize
  \begin{equation} \label{kmeans-sam}
 	\Phi_n(a_1,\ldots,a_k) =  \Phi({\mathbb P}_n;a_1,\ldots,a_k) = \frac{1}{n}\sum_{j=1}^n\left(
 	\min_{i=1,...,k} \|X_j -a_i\|_{\mathcal{B}}^2 \right).
 \end{equation}
 The values minimizing \eqref{kmeans-sam} are called the \textit{empirical $k$-means}. The empirical WCSS is then
 \begin{equation}\label{empirical-WCSS}
  W_n(k)= W({\mathbb P}_n;k)= \inf_{\{a_1,\dots,a_k\}\subset\mathcal{B}} \Phi_n(a_1,\ldots,a_k).
  \end{equation}

The $k$-means procedure plays a central role in localization problems in operations research.
It can be motivated in terms of the so-called \textit{facility location problem}: determine the optimal placement, $\mu_1,\dots,\mu_k$, for $k$ facilities in such a way that the average (square) distance from a random individual to the closest facility is minimal. In statistics and machine learning, $k$-means techniques are generally used in clustering (see, e.g., \cite{jai10} for a survey), where the aim is partitioning the space in a \textit{Voronoi tesellation} of cells associated with the $k$-means. Hence, we divide the sampling space (or the sample) into a partition of $k$ clusters. The $i$-th cluster is constituted by all elements whose closest center is $\mu_i$ ($i=1,\dots,k$).

In addition to the usual and multiple applications of this technique (unsupervised classification, taxonomy, image analysis, information
retrieval, market segmentation, computer vision, astronomy, etc.), $k$-means clustering is also used in \textit{quantization}: constrain a possibly infinite large set to a small collection of values, the $k$-centers. In particular, the $k$-means method is extensively used by practitioners as a powerful tool for summarizing data in cluster prototypes. The group centers provide a basic representation of the data that is usually very informative. The information contained in the principal points is extremely useful as a descriptive tool as it often allows understanding the underlying structure of the data and identifying prominent features.

The $k$-means methodology has important advantages: it is completely general, applicable to data in normed or metric spaces; easy to interpret and understand; with  population/sampling versions perfectly defined; several efficient heuristic algorithms are available; and with a sound asymptotic theory behind it. In spite of these appealing aspects, this popular technique also entails some challenges. First, the effective calculation of the sample $k$-means is a formidable computational (NP-hard) problem: algorithms have to cope with a non-convex optimization problem in a possibly high or even infinite-dimensional space. Moreover, the usual algorithms do not guarantee to reach a global optimum, but rather a local one; see \cite{Ikotun-2023} and the references therein
for a recent overview of various clustering algorithms. Choosing of a good value for $k$ still receives considerable attention and constitutes an area of active research: it is difficult to give a clear solution on the choice of the ``best'' number of groups, $k$, whatever method is used. Finally, as we discuss throughout this paper, the $k$-means problem does not necessarily have a single solution in the population. The lack of uniqueness might lead to important stability problems of the algorithms, as well as more subtle issues that can affect the understanding and interpretation of the results.

The uniqueness of the $k$-means set is imposed as a condition for obtaining results in many works on this subject. Still, it is not difficult to find simple examples where this uniqueness prerequisite is violated. In Section \ref{sec:non-unique} we present some models to illustrate that there are in fact at least two distinct cases of $k$-means non-uniqueness that should be differentiated.  
Section \ref{sec:main} brings together the main theoretical results. We show a general consistency theorem that does not assume uniqueness using the Hausdorff metric and the Gromov-Hausdorff distance. We also obtain the asymptotic distribution of the empirical WCSS in \eqref{empirical-WCSS}. This asymptotic result relies on the class of functions generating the $k$-means problem verifying the Donsker property. At the end of Section \ref{sec:main} we provide some technical results that guarantee this property. 
In Section \ref{sec:test} we prove that the uniqueness of the population $k$-means set is equivalent to the asymptotic normality of the empirical WCSS, as well as other useful characterizations. As an application, we propose a bootstrap test for the uniqueness of the set of $k$-means. Finally, in Section~\ref{sec:emp} we collect some numerical simulations to check the performance of our proposal in practice.


\section{Patterns of (non-)uniqueness: a zoo of examples}\label{sec:non-unique}

The purpose of this section is to gain some insight into the phenomenon of non-uniqueness in $k$-means. We discuss two different patterns of non-uniqueness and provide a collection of illustrative examples. These examples are used in Section \ref{sec:emp} for numerical simulations.


\subsection{Basic definitions and notation}

We start with some definitions that we use throughout this work. We denote by $\mathcal{S}_{\P}(k)$ the collection of all the $k$-means sets of the probability measure $\P$, that is,
\begin{equation}\label{k-means-set}
 \mathcal{S}_{\P}(k) =  \{ K=\{\mu_1,\dots,\mu_k\}\subset{\mathcal{B}} : \text{ $K$ is a set of $k$-means for }  \P  \}.
\end{equation}
We say that $\P$ satisfies the \textit{$k$-means uniqueness property}, UP($k$), or $\P \in \text{UP}(k)$ in short, if the set $\mathcal{S}_{\P}(k)$ in \eqref{k-means-set} has cardinality one. Therefore, if $\P \in \text{UP}(k)$, there is a single set of $k$-means attaining the infimum of the population WCSS in \eqref{WCSS-population}. Historically, this property has been an initial (and essential) assumption in all significant results that provide theoretical support for this methodology; see the seminal paper by \cite{pol81} on the consistency of the procedure and \cite{pol82} where the asymptotic normality of the empirical $k$-means is established. This hypothesis has been maintained in all subsequent works on this subject; see for example \cite{cue88} and \cite{lem03}.

Further, as pointed out by \cite{Garcia-Escudero-et-al-1999}, checking the $k$-means uniqueness is a difficult task in practice. There are only a few references related to the ``uniqueness of the principal points''; see \cite{Li-Flury-1995}, \cite{Tarpey-1994}, \cite{Trushkin-1982} and \cite{Zoppe-1997}. These contributions mostly deal with very particular cases or univariate distributions which are not too relevant in clustering analysis. In general, it is not easy to determine analytically whether a given multidimensional distribution verifies the uniqueness property UP($k$) or not.

In the machine learning literature, it is commonly accepted that the lack of uniqueness is equivalent to $k$-means algorithms having instability problems. In short, it is accepted that the existence of a unique minimizer amounts to the stability of the $k$-means clustering; see \cite{Ben-David-2006}, \cite{Ben-David}, \cite{Caponnetto-Rakhlin-2006} and the overview by \cite{Von Luxburg}. This is of some practical significance since stability might be used for choosing the number of clusters as $k$ can be selected as the value that provides the most stable results; see \cite{Caponnetto-Rakhlin-2006}. However, as we point out in Section \ref{sec:emp}, this equivalence between uniqueness and stability is not entirely accurate because there are situations of non-uniqueness in which the algorithms are shown to be stable. This occurs when the different sets of $k$-means are ``separated'' from each other. This idea is elaborated in what follows.


We first note that the cases of non-uniqueness usually appear in two different modalities. To discuss this we use the classical Hausdorff distance. Given two non-empty compact sets $A,C\subset {\mathcal B}$, the \textit{Hausdorff distance} between these two sets is
\begin{equation}\label{Haus}
d_H(A,C)=\inf\{\varepsilon>0: A\subset C^\varepsilon \mbox { and } C\subset A^\varepsilon\},
\end{equation}
where, for any set $B\subset {\mathcal B}$, $B^\varepsilon=\bigcup_{x\in B}\{z: \Vert x-z\Vert\leq \varepsilon\}$ is the $\varepsilon$-dilation of $B$. It is well-known (see, e.g., \cite{mun74}) that if $\mathcal B$ is a separable Banach space, then
${\mathcal H}({\mathcal B})= \left\{ B\subset \mathcal{B} : \emptyset \ne B \text{ compact} \right\}
$
is a complete separable metric space, when endowed with the Hausdorff metric $d_H$ in \eqref{Haus}.

We distinguish two different types of \textit{multiplicity patterns}, that is, different situations in which there is no uniqueness in the $k$-means problem.

\begin{itemize}
    \item [--] \textsl{Continuous non-uniqueness}, $\textrm{CNU}(k)$: For each $\varepsilon>0$ and each set of $k$-means, $K_1\in \mathcal{S}_\P(k)$, there exists $K_2\in \mathcal{S}_\P(k)$ such that $K_1\ne K_2$ and $d_H(K_1,K_2)<\varepsilon$. Here, every $k$-means set is an accumulation point (with respect to the Hausdorff metric) of the set $\mathcal{S}_\P(k)$.

    \item [--] \textsl{Discrete non-uniqueness}, DNU($k$): The set of $k$-means  $\mathcal{S}_\P(k)$ has cardinality greater than one and there exists $\varepsilon_0>0$ such that $d_H(K_1,K_2)\ge \varepsilon_0$, for each pair $(K_1,K_2)$ of different sets in $\mathcal{S}_\P(k)$.
\end{itemize}

Whenever $\P\in \text{CNU}(k)$, there are different sets of $k$-means arbitrarily close to each other. In practice, this is the worst possible situation because it effectively leads to $k$-means algorithms needing many iterations and a large amount of time and computational cost to stop. This is the case that is often identified by the machine learning community with the instability of the $k$-means algorithms. In contrast, in the DNU($k$) case, there are multiple but isolated $k$-means sets. In practice, when sampling from $\P\in  \text{DNU}(k)$, the algorithms usually approach one of the $k$-means sets and converge without instability problems. This difference between the two types of non-uniqueness patterns has been corroborated by numerical simulations in Section \ref{sec:emp}: we have observed that the average run times for $\text{CNU}(k)$ models are significantly higher than for the case DNU($k$).

In the remainder of this section we describe several models of the different cases of non-uniqueness. With these examples we also show the complex casuistry that can arise in $k$-means, even in the case of uniqueness. In particular, we emphasize that uniqueness does not necessarily imply a good choice of the parameter $k$.
These examples are revisited in Section \ref{sec:emp} to illustrate the theoretical results in Section \ref{sec:main} and evaluate the practical performance of the uniqueness test proposed in Section \ref{sec:test}.

\subsection{Examples of different types of non-uniqueness}


In all the considered models we use the following notation: C$i$k$j$ stands for a mixture of $i$ distributions (number of components) and we decide to look for $k=j$ centers. As there may be several examples of the same case (equal number of components and value of $k$), in some instances we will add some extra letter to specify the type of model we are using; U (uniform), T (triangle), etc. In this section we only consider distributions in $\R$ or $\R^2$ to facilitate the visualization of the probability densities and the associated $k$-means sets. Therefore, we do not use here any dimension identifier. However,  most of these models will be extended to higher dimensions in Section \ref{sec:emp}.

We begin with a very simple one-dimensional (toy) example in which the difficulty of the $k$-means problem can be perceived. 

\begin{example}[Model U$r$C3k2]\label{Example-uniform} {\rm For $0 \le r \le 1/2$, let $\P_r$ be a probability measure of a mixture of three uniform distributions in $\R$ with density
\begin{equation}\label{density-uniform}
  f_r(x)= \frac{1}{6 r}\left(  1_{B(-1,r)} (x)  +  1_{B(0,r)}(x)  +1_{B(1,r)} (x) \right), \quad 0<r\le \frac{1}{2}, \quad x\in \R,
\end{equation}
where $1_A$ stands for the indicator of the set $A$ and $B(x,r)$ is an open ball centered at $x$ of radius $r$. The case $r=0$ corresponds to the discrete uniform distribution on $\{-1,0,1\}$ and $r=1/2$ to the uniform distribution on the interval $(0,1)$.

We select $k=2$. Even in this elementary example it is not straightforward to determine the sets of 2-means. After some calculations, we obtain that
\begin{equation*}
\mathcal{S}_{\P_r}(2)=
  \begin{cases}
    \{ \{ -1,1/2 \}, \{-1/2,1\}  \}, & \mbox{if } 0\le r < 3\sqrt{2}-4, \\
    \{ \{ -1,1/2 \}, \{-1/2,1\}, \{-1/\sqrt{2}, 1/\sqrt{2} \} \}, & \mbox{if } r= 3\sqrt{2}-4, \\
     \{ \{ -(4+r)/6 ,  (4+r)/6 \}  \}  ,  & \mbox{if } 3\sqrt{2}-4< r \le 1/2.
  \end{cases}
\end{equation*}

The population WCSS can also be calculated and we have that
\begin{equation*}
W(\P_r;2)=
  \begin{cases}
   {\displaystyle\frac{2r^2+
   1}{6}}, & \mbox{if } 0\le r \le 3\sqrt{2}-4, \\[3 mm]
     \displaystyle\frac{11 r^2 -8 r +8}{36} ,  & \mbox{if } 3\sqrt{2}-4< r \le 1/2.
  \end{cases}
\end{equation*}

\begin{figure}[h]\label{fig:UrC3k2-population}
\centering
     \includegraphics[width=0.90\textwidth]{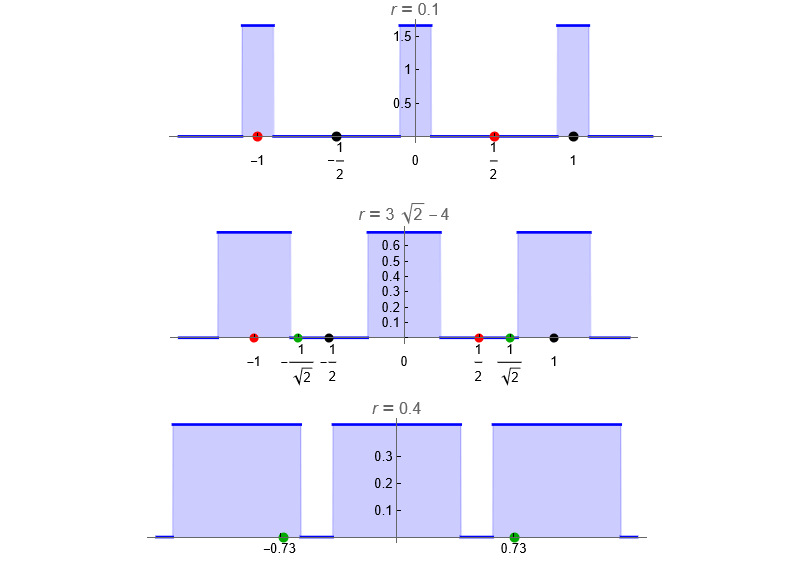}
    \caption{Model U$r$C3K2. Densities functions in \eqref{density-uniform} with the  sets of $2$-means (in black, red and green) for $r=0.1$ (upper panel),  $r=3\sqrt{2}-4$ (middle panel) and $r=0.4$ (lower panel).}
\end{figure}

Figure \ref{fig:UrC3k2-population} shows three graphs of $f_r$ in \eqref{density-uniform}, along with their sets of 2-means. Observe that $\P_r\in \textrm{DNU}(2)$ if $r\in[0, 3\sqrt{2}-4]$.
Just at $r=3\sqrt{2}-4$ there is a ``phase change'' from non-uniqueness to uniqueness and $\P_r\in \textrm{UP}(2)$ for $r\in(3\sqrt{2}-4, 1/2]$. We note that uniqueness does not entail a good choice of $k$ for $r\in(3\sqrt{2}-4, 1/2)$. Moreover, for $r\in[0, 3\sqrt{2}-4)$ (non-uniqueness), we can go from one set of 2-means to another one using a symmetry with respect to the origin.
}
\end{example}

\begin{example}[Model C1k2] \label{Example-C1k2} {\rm
One of the simplest examples of the CNU case is obtained by choosing a value of $k$ greater than 1 in a population with a radial density function. We consider $\operatorname{P}\sim\mathcal{N}\left((0,0),\mathbb{I}_{2}\right)$,
distributed as a standard bivariate normal, where $\mathbb{I}_{2}$ is the identity matrix of dimension 2. When we select $k=2$, due to the circular symmetry of the normal density, infinite sets of $k$-means appear. It is not difficult to check that these sets are formed by any two diametrically opposite points on a circumference centered at the origin of radius $\sqrt{2/\pi}\approx 0.16$; see Figure \ref{fig:C1k2-population}.

There is an symmetry pattern in $\mathcal{S}_\P(2)$ (the collection of $2$-means sets): any set of 2-means can be brought to another one by a rotation around the origin. The same happens
for other values of $k$ greater than 2.

\begin{figure}[h]\label{fig:C1k2-population}
\centering
    \begin{tabular}{cc}
     \includegraphics[width=0.62\textwidth]{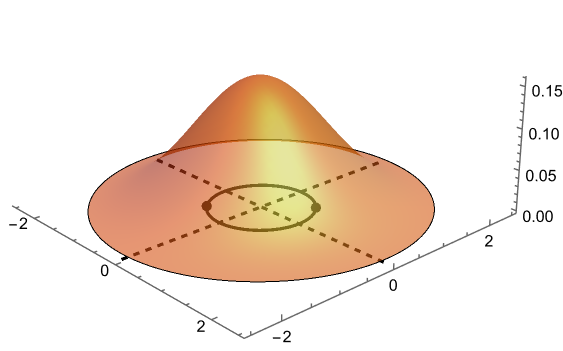} &
    \includegraphics[width=0.32\textwidth]{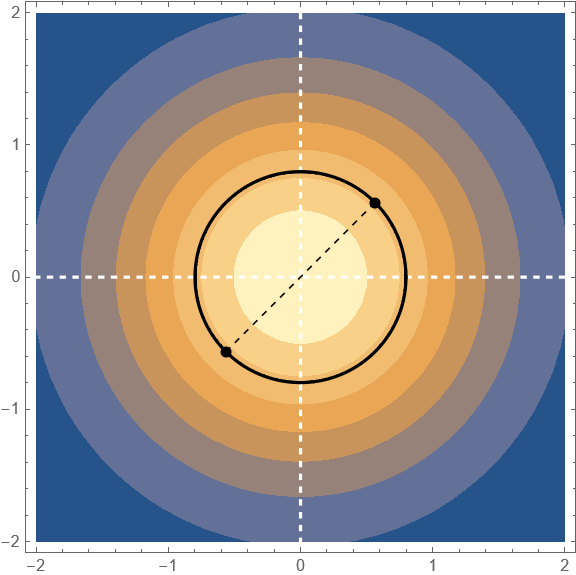} \\
    (a) & (b)
    \end{tabular}
    \caption{Model C1k2. (a) Density plot and (b) contour plot. The two points in black are one of the sets of $2$-means located on a circumference centered at the origin of radius $\sqrt{2/\pi}$.
    }
    \label{fig:C1k2-population}
\end{figure}
}
\end{example}

\begin{example}[Model C2k3]\label{ExampleC2k3} {\rm
Another example of continuous non-uniqueness is obtained when selecting $k=3$ with a suitable mixture of two bivariate normal distributions. Specifically, we consider
\begin{equation}\label{model-C2k3D2}
\operatorname{P}\sim\frac{1}{2}\,\mathcal{N}\left((-1,0), \frac{\mathbb{I}_{2}}{25}\right)+\frac{1}{2}\,\mathcal{N}\left((1,0),\,\frac{\mathbb{I}_{2}}{25}\right).
\end{equation}

To facilitate the computation of the 3-means we assume that the normal distributions are truncated from 5 deviations. In this way, the tails of one component do not affect the other and infinite sets of $3$-means appear of the case $\text{CNU}(3)$. It can be seen that the sets are formed by two diametrically opposite points on a circumference centered on the mean of one of the normal distributions of the mixture of radius $\sqrt{2/\pi}$ together with the center of the other normal. A graphical representation of this situation is presented in Figure \ref{fig:C2k3-population}. Note that there is an isometric transformation (a rotation plus a reflection) that takes one of these $3$-means sets into another.

\begin{figure}[h]\label{fig:C1k2-population}
\centering
    \begin{tabular}{c}
     \includegraphics[width=0.70\textwidth]{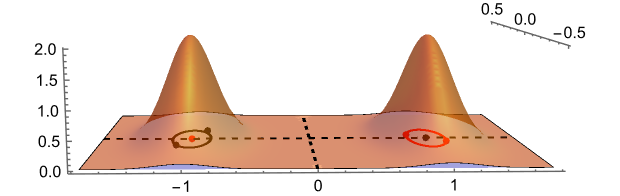} \\
     (a) \\[3 mm]
    \includegraphics[width=0.62\textwidth]{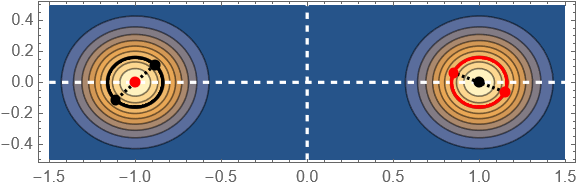} \\
    (b)
    \end{tabular}
    \caption{Model C2k3. (a) Density plot and (b) contour plot. The two sets of three points in black and red are two of the infinite sets of $3$-means of the case $\text{CNU}(3)$.}
    \label{fig:C2k3-population}
\end{figure}

}
\end{example}

\begin{example}[Model TC3k2, Gaussian triangle]\label{Example-TC3k2} {\rm We consider an example of discrete non-uniqueness with Gaussian variables. We set $k=2$ in a mixture of three bivariate normal distributions (with equal weights). Let
$$\operatorname{P}\sim\frac{1}{3}\,\mathcal{N}\left((1,0),\frac{\mathbb{I}_{2}}{25}\right)
+
\frac{1}{3}\,\mathcal{N}\left(\left(-\frac{1}{2},\frac{\sqrt{3}}{2}\right), \frac{\mathbb{I}_{2}}{25}\right)
+
\frac{1}{3}\,\mathcal{N}\left( \left(-\frac{1}{2},-\frac{\sqrt{3}}{2}\right),\,\frac{\mathbb{I}_{2}}{25}\right).$$

Note that the centers of the 3 components are the cube roots of unity, so that the set of centers is invariant by rotations of $\pi/3$ and $2\pi/3$ radians. Therefore, by choosing $k=2$ we obtain 3 distinct sets of 2-means; see Figure \ref{fig:TC3k2-population}. 

\begin{figure}[h]    \label{fig:TC3k2-population}
\centering
    \begin{tabular}{c}
     \includegraphics[width=0.8\textwidth]{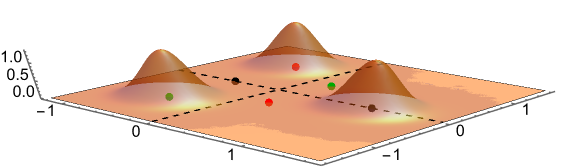} \\
     (a) \\[3 mm]
    \includegraphics[width=0.35\textwidth]{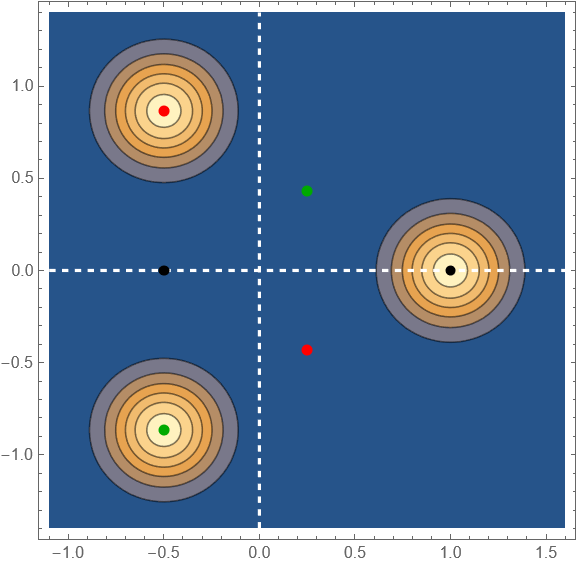} \\
    (b)
    \end{tabular}
    \caption{Model TC3k2. (a) Density plot and (b) contour plot. The three sets of 2 points in black, red and green are the sets of $2$-means of the case $\text{DNU}(2)$.}

\end{figure}

}
\end{example}

The uniqueness property can easily fail in higher dimensions. In fact, it seems easier for the symmetries observed in the above examples to occur when the dimension of the underlying space is higher. As a consequence, $k$-means multiplicity could appear more frequently in those situations where geometric and visual intuition is of little help. Simple models not satisfying the UP($k$) can also be found in the infinite-dimensional setting, i.e., in the context of functional data analysis. For instance, let us consider $(t,B_1(t),B_2(t))$, where $(B_1(t),B_2(t))$ is a vector of independent and standard Brownian motions and $t\in [0,1]$ (time variable), and let us fix $k\ge 2$. In this case, there is a symmetry around the time-axis and infinite $k$-means sets appear of the case $\text{CNU}(k)$.

It is important to note that in the above examples the multiplicity in the $k$-means appears when: (1) the distribution of the population has a certain symmetry; and (2) the value of $k$ has not been suitably chosen. In practice, it is very relevant to detect either of these two circumstances. 

\subsection{Examples of different types of uniqueness}

The following examples include different situations in which the uniqueness property is fulfilled.

\begin{example}[Models C2k2] \label{Example-C2k2} {\rm We consider the mixture of two Gaussian distributions and we correctly select the value $k=2$. Therefore, this is a simple situation in which we have uniqueness with a good choice of the parameter $k$. We consider 3 homocedastic models that differ in the (Mahalanobis) distance between the vectors of means.

\begin{itemize}
  \item[--] Model C2k2-1: The measure $\P$ is given in \eqref{model-C2k3D2} (see Example \ref{ExampleC2k3}).
The two populations in the mixture are well-separated with density plotted in Figure \ref{fig:C1k2-population}. Some calculations show that the set of $2$-means is $\mathcal{S}_{\P}(2)=\{ (-c_1,0 ),(c_1,0)  \}$, where
\begin{equation}\label{c-case-c2k2}
c_1= \frac{e^{-25/2}}{15} \sqrt{\frac{2}{\pi}} + 2\Phi(5)-1,
\end{equation}
where $\Phi(\cdot)$ is the cdf of a standard normal distribution in $\R$. Observe that the set of 2-means is not exactly the means of the two components of the mixture. Actually, the value $c_1$ in \eqref{c-case-c2k2} is slightly greater than 1. This is due to a kind of \textit{repulsion effect} between the centers because the tail of every component affects the center of the other one.

\item[--] Model C2k2-2: We consider a slight modification of the previous model in which the centers are at distance 3. In this case, $\P$ is the mixture 
given by
\begin{equation}\label{model-C2k3D2}
\operatorname{P}\sim\frac{1}{2}\,\mathcal{N}\left((-3/2,0), \mathbb{I}_{2}\right)+\frac{1}{2}\,\mathcal{N}\left((3/2,0),\,{\mathbb{I}_{2}}\right).
\end{equation}
Here, there is a single set of $2$-means given by $\{(-c_2,0),(c_2,0)\}$, where
\begin{equation*}
c_2= e^{-9/8}\sqrt{\frac{2}{\pi}} + \frac{3}{2}(2\Phi(3/2)-1)\approx 1.56.
\end{equation*}
The density of this distribution and the set of $2$-means are shown in Figure \ref{fig: C2k2-Mofified}.

 \item[--] Model C2k2-3:  Let  $\P$ be the mixture of two Gaussian distributions given by
\begin{equation}\label{model-C2k3D2}
\operatorname{P}\sim\frac{1}{2}\,\mathcal{N}\left((-1,0), \mathbb{I}_{2}\right)+\frac{1}{2}\,\mathcal{N}\left((1,0),\,{\mathbb{I}_{2}}\right).
\end{equation}
In this example, the centers of the two populations are closer to each other and the joint density has the appearance of a single mountain with a unique maximum at the point $(0,0)$. Still, there is a single set of $2$-means given by $\{(-c_3,0),(c_3,0)\}$, where
\begin{equation*}
c_3= \sqrt{\frac{2}{e \pi}} + 2\Phi(1)-1\approx 1.17.
\end{equation*}
We expect this uniqueness case to be more difficult to identify by numerical simulations than the previous models C2k2-1 and  C2k2-2 since we are closer to the multiplicity of the case C1k2 in Example \ref{Example-C1k2}. We also observe that the repulsion effect of the 2-means is more evident when the centers are closer together; see Figure \ref{fig: C2k2-Mofified}

\end{itemize}

\begin{figure}[h]\label{fig: C2k2-Mofified}
\centering
    \begin{tabular}{c c}
     \includegraphics[width=0.48\textwidth]{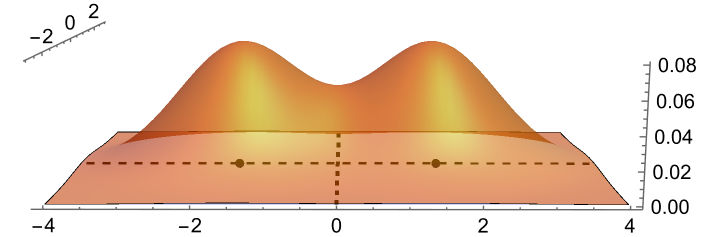} &  \includegraphics[width=0.48\textwidth]{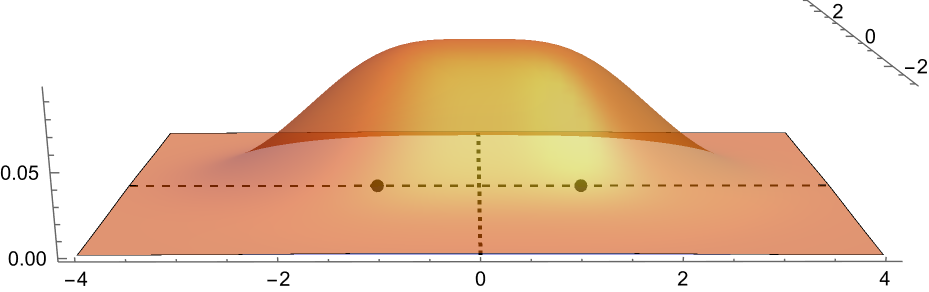}  \\
     (a) & (c) \\[3 mm]
    \includegraphics[width=0.45\textwidth]{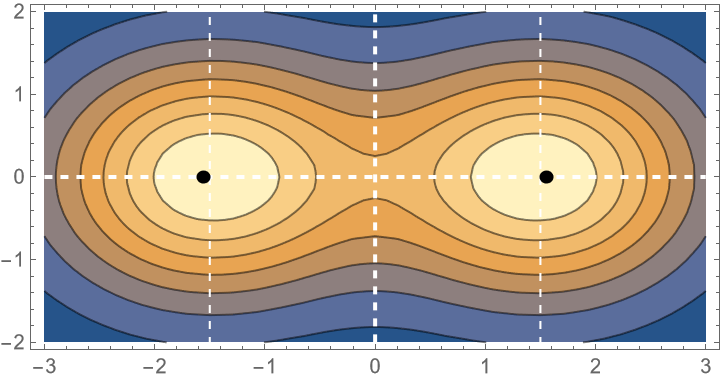} &  \includegraphics[width=0.45\textwidth]{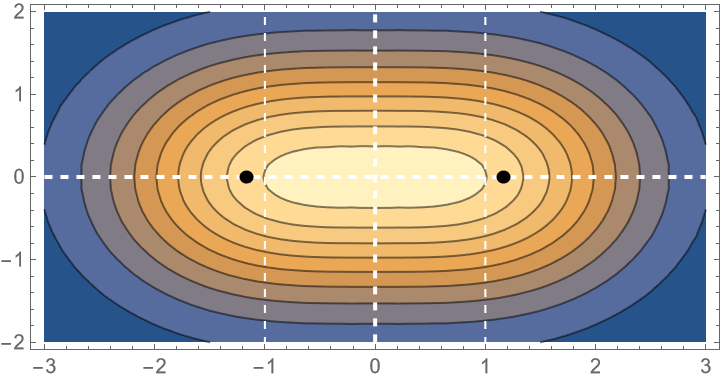} \\
    (b) & (d)
    \end{tabular}
    \caption{Models C2k2-2 (left panels) and C2k2-3 (right panels). Density plots in (a) and (c) and  contour plots in (b) and (d). The points in black are the (unique) sets of $2$-means.}
\end{figure}

}
\end{example}

\begin{example}[Model C3k3]\label{Example-C3k3} {\rm Let $X$ be a mixture of three bivariate normal distributions (with equal weights) and probability measure
\begin{equation}\label{C3k3-density}
\operatorname{P}\sim\frac{1}{3}\,\mathcal{N}\left((-1,0),\frac{\mathbb{I}_{2}}{25}\right)
+
\frac{1}{3}\,\mathcal{N}\left((0,0), \frac{\mathbb{I}_{2}}{25}\right)
+
\frac{1}{3}\,\mathcal{N}\left((1,0),\,\frac{\mathbb{I}_{2}}{25}\right).
\end{equation}
The density of this variable is represented in Figure \ref{fig:C3k2-population}.
For the value $k=3$, we obtain that $\P\in\text{UP(3)}$. The sets of $3$-means can be computed analytically and it is given by $\{ (-c,0), (0,0), (c,0)  \}$ with a value of $c$ slightly greater than 1 (due again to the repulsion effect mentioned above). The exact value of $c$ is omitted for the sake of brevity. This is, as Example \ref{Example-C2k2}, a case of uniqueness with a good choice of $k$.
}
\end{example}

\begin{example}[Model C3k2]\label{Example.last} {\rm
We consider an example of uniqueness in which the choice of parameter $k$ is wrong. Specifically, we select $k=2$ in the population given in \eqref{C3k3-density} (as in Example \ref{Example-C3k3}). 
After some computations, we obtain the value of the unique set of 2-means given by $\{(-c,0),(c,0)\}$, with
\begin{equation*}
c= \frac{1}{15} \left(10 (2\Phi(5)-1) + \sqrt{\frac{2}{\pi }}      (1+2 e^{-25/2}) \right)\approx 0.72.
\end{equation*}
Therefore, $\operatorname{P}\in \text{UP}(2)$. See Figure \ref{fig:C3k2-population} for a graphical representation of this example. As we will see in Section \ref{sec:emp}, this typology is more difficult to detect by numerical simulations than other situations of uniqueness in which $k$ is properly selected.

\begin{figure}[h]\label{fig:C3k2-population}
\centering
    \begin{tabular}{c}
     \includegraphics[width=0.62\textwidth]{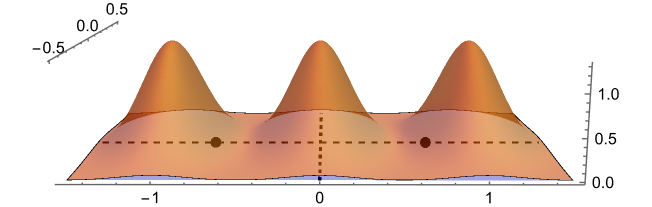} \\
     (a) \\[3 mm]
    \includegraphics[width=0.62\textwidth]{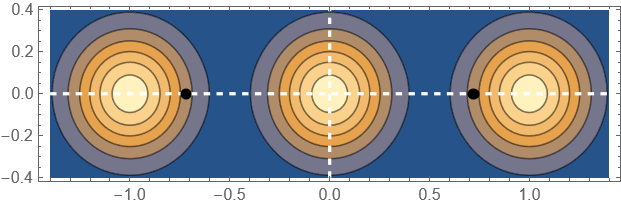} \\
    (b)
    \end{tabular}
    \caption{Model C3k2. (a) Density plot and (b) contour plot. The set of two points in black is the unique set of $2$-means, case $\text{UP}(2)$.}
    \label{fig:C3k2-population}
\end{figure}
}
\end{example}




\section{Main theoretical results}\label{sec:main}

In this section we show the consistency of the empirical $k$-means in terms of the Gromov-Hausdorff metric even under non-uniqueness. We also derive the asymptotic distribution of the empirical WCSS. Finally, we provide various sufficient conditions for the class of functions defining the $k$-means risk minimization problem to be Donsker. This last result is important from a theoretical point of view because it is the basic assumption to obtain the characterization of the uniqueness property in Section \ref{sec:test}.

\subsection{Strong consistency without uniqueness}\label{subsec:Strong}

As we have shown in the previous section, some simple models might not satisfy the uniqueness property UP($k$). Particularly in the $\text{CNU}(k)$ case, the standard $k$-means algorithms could show a remarkable instability, providing different outputs, depending on the initial conditions or failing to fulfil the standard stopping criteria. However, if the sample size $n$ grows to infinity and we obtain a set of empirical $k$-means for each $n$, one might ask about the ``limit behaviour'' of such sequence of sets. The following Theorem \ref{th_cluster_consistency} provides an answer. It is an elaboration from  some previous results by \cite{cue88} and \cite[Th. 3.1]{lem03}. Essentially, the result establishes that the set of empirical $k$-means always approaches some population $k$-mean set and, reciprocally, every population $k$-mean set is a limit of a (sub)sequence of empirical $k$-means sets. Here, the results of convergence for sets of $k$-means is stated in terms of the Hausdorff and the Gromov-Hausdorff metric, both defined on the space ${\mathcal H}({\mathcal B})$ of compact non-empty sets introduced above.

The \textit{Gromov-Hausdorff metric} (see \cite{bur01} for details) is defined by
\begin{equation}\label{GH}
	d_{GH}(A,C)=\inf\{d_H(T(A),S(C))\},\quad A,C\in {\mathcal H}({\mathcal B}),
\end{equation}
where $d_H$ is a the Hausdorff metric in \eqref{Haus} and the infimum ranges over all possible choices of $T$, $S$ and ${\mathcal M}$, with $T:A\rightarrow {\mathcal M}$ and $S: C\rightarrow {\mathcal M}$ being isometric embeddings and ${\mathcal M}$ is a compact metric space.

We use the following assumptions:

\begin{description}

\labitem{(Geo)}{itm:Geo} \textit{Geometric assumption.} ${\mathcal B}$ is a separable and uniformly convex Banach space.

\labitem{(Int)}{itm:Int} \textit{Integrability assumption.} $X$ satisfies $\E\Vert X\Vert_{\mathcal{B}}^2<\infty$ (strong second moment).

\labitem{(Sym)}{itm:Sym} \textit{Symmetry assumption.} For any $K_1,K_2\in {\mathcal S}_{\P}(k)$ there is an isometry $T:{\mathcal B}\rightarrow {\mathcal B}$ such that $T(K_2)=K_1$.

\end{description}

Assumption  \ref{itm:Geo} requires that $\mathcal{B}$ is uniformly convex or uniformly rotund; see \cite{cla36}. This is fulfilled for Hilbert spaces as well as for $\mathrm{L}^p$ spaces, with $1<p<\infty$. This condition is imposed by \cite{cue88} to derive their consistency results. This assumption entails the less restrictive, but more technical, hypothesis in \cite[Th. 3.1]{lem03}. The requirement \ref{itm:Int} is necessary for the functional \eqref{kmeans-pop} to be well-defined. Assumption \ref{itm:Sym} might seem somewhat restrictive. However, it holds in almost all of the examples in Section \ref{sec:non-unique} where the sets of $k$-means present self-similarity patterns easy to formalize in terms of isometries (translations, rotations, reflections, \dots).


\begin{theorem}\label{th_cluster_consistency}
Let $X$ be a random element taking values in ${\mathcal B}$  with probability distribution $\P$ fulfilling \ref{itm:Geo} and \ref{itm:Int}. We consider the sequence $(K_n)$ (for $n\in\mathbb{N}$), where $K_n\in \mathcal{S}_{\mathbb{P}_n}(k)$ is a $k$-mean set of the empirical measure $\mathbb{P}_n$. We have that:

 \begin{itemize}
     \item [(i)] Any $d_H$-adherent point of the sequence $(K_{n})$ belongs to $\mathcal{S}_{\P}(k)$ a.s.
\item [(ii)] If additionally \ref{itm:Sym} holds, for any $K_0\in\mathcal{S}_{\P}(k)$, there is a subsequence of $(K_n)$, say $(K_{n_j})$ (for $j\in\mathbb{N}$), such that $d_{GH}(K_{n_j},K_0)\to 0$, almost surely, as $j\to\infty$.
 \end{itemize}
\end{theorem}

\begin{proof}

Part (i) essentially follows from \citet[Th. 3.1]{lem03}. Such result is very general but quite technical. Therefore, we check all its requirements below.
To begin with, we consider the ordinary weak topology in $\mathcal{B}$ as $\tau$. Then, we first prove that the key assumption \cite[Assumption \textbf{B}, p. 29]{lem03} is fulfilled for this choice of $\tau$. This condition  imposes that ``every closed ball of ${\mathcal B}$ is sequentially $\tau$-compact''. This is guaranteed by \ref{itm:Geo}. Indeed, from Milman-Pettis theorem, every uniformly convex space is reflexive (we refer to \citet[Ch. 3]{bre10} for the functional analysis results used in this proof) and hence
the weak and the weak$^*$ topology coincide. From Banach-Alaoglu theorem, the unit ball in $\mathcal{B}$ is compact in the weak$^*$ topology. In addition, as $\mathcal{B}$ is separable, the sequential version of Banach-Alaoglu theorem holds, so that every closed ball in $\mathcal{B}$ is sequentially compact in the  weak$^*$ topology, and also in the weak one (that we denote $\tau$) by the reflexive property. So we conclude that \textbf{B} assumption in \citet[Th. 3.1]{lem03} is fulfilled.

Also, Assumption (2) in Lember's result holds here. This is the so-called Kadec-Klee property that, according to \citet[p. 29]{lem03} amounts to the following Radon-Riesz property, when $\tau$ is the weak topology: if $x_n$ converges weakly to $x$ in ${\mathcal B}$ and $\|x_n\|_{\mathcal{B}}\to \|x\|_{\mathcal{B}}$, then $\|x_n-x\|_{\mathcal{B}}\to 0$. But this property holds in uniformly convex Banach spaces \cite[Prop. 3.32]{bre10}.

As a conclusion, we can apply \cite[Th. 3.1]{lem03} to conclude that every subsequence of the empirical $k$-means $( K_n)$ has a further subsequence converging (almost surely) in the Hausdorff metric to some  set $K_0\in \mathcal{S}_{\P}(k)$. Observe that if $(K_{n_j})$ is a subsequence of the empirical $k$-means converging almost surely to $K_0\notin \mathcal{S}_{\P}(k)$, all subsequences of $(K_{n_j})$ should converge almost surely to $K_0$, which contradicts \cite[Th. 3.1]{lem03}.

To prove (ii), we consider $K_0\in\mathcal{S}_{\P}(k)$. Let $(K_n)$ be a sequence with $K_n\in \mathcal{S}_{\mathbb{P}_n}(k)$. Using (i), there exists a subsequence, say $(K_{n_j})$, and $K\in{\mathcal S}_\P(k)$ such that $d_H(K_{n_j}, K)\to 0$, almost surely, as $j\to\infty$. Now, using \ref{itm:Sym}, let $T$ be an isometry on ${\mathcal B}$ such that $T(K)=K_0$. Let us denote by ${\mathcal X}_{n_j}$ the subsequence of data sets corresponding to $K_{n_j}$. Since $T$ is an isometry, a sequence of $k$-means sets corresponding to the sequence $T({\mathcal X}_{n_j})$ is $T(K_{n_j})$. As $d_H(K_{n_j},K)\to 0$ a.s., we also obtain that
\begin{equation}\label{star}
d_H(T(K_{n_j}),K_0)=d_H(T(K_{n_j}),T(K))\stackrel{(*)}{=}d_H(K_{n_j},K)\rightarrow 0, \ \mbox{a.s.}
\end{equation}
To see ($\ast$), recall that, given $A, C\in {\mathcal H}({\mathcal B})$, the Hausdorff distance in \eqref{Haus} between them can be alternatively expressed as
\begin{equation}\label{Haus2}
	d_H(A,C)=\max\left\{\sup_{a\in A}d(a,C),\, \sup_{c\in C}d(c,A)\right\},
\end{equation}
where $d(a,C)=\inf_{c\in C}\|a-c\|_{\mathcal{B}}$ and similarly for $d(c,A)$. As a consequence of \eqref{Haus2}, Hausdorff metric remains invariant when the same isometry is applied to both sets. So, equality ($\ast$) in \eqref{star} holds.

Finally, from the definition \eqref{GH} of $d_{GH}$ and \eqref{star} we obtain that
$$
d_{GH}(K_{n_j},K_0)\leq d_H(T(K_{n_j}),K_0)\to 0,\quad \text{a.s. as } j\to\infty.
$$
This concludes the proof.

\end{proof}



Theorem \ref{th_cluster_consistency} can be simplified if we could guarantee the (almost sure) uniqueness of the sequence of empirical $k$-means.
Intuitively, the uniqueness, with probability one, of the $k$-means of the sample seems quite natural whenever the distance $\Vert X-a\Vert_{\mathcal{B}}$ ($a\in \mathcal{B}$) has an absolutely continuous distribution. However, to the best of our knowledge, this ``empirical uniqueness'' does not seem simple to prove. It is essentially a matter of establishing the uniqueness of a sample sequence defined in terms of the `$\argmin$' of a suitable functional. The interesting paper by \cite{cox20} deals with the uniqueness of $\argmin$-type statistics with a especial focus of M-statistics and maximum likelihood estimators. The methodology in this paper includes differentiability assumptions which make little sense in the $k$-means framework. This suggests us that the full study of affordable sufficient conditions for the uniqueness of the empirical $k$-means is an interesting question far beyond the scope of the present work.

\subsection{Asymptotic behaviour of the empirical WCSS}

In this section $\mathcal{B}$ is a Banach space, without the constraint \ref{itm:Geo} used in Section~\ref{subsec:Strong}. The $k$-means problem, as stated in \eqref{kmeans-pop} and \eqref{kmeans-sam}, can be viewed as a risk minimization problem  over an appropriate class of functions. We follow this approach here, within the context of empirical processes theory, to derive the asymptotic distribution of the empirical WCSS ($W_n(k)$ in \eqref{empirical-WCSS}), and subsequently address the question of the uniqueness in $k$-means in Section~\ref{sec:test}.
In what follows we use the notation $\operatorname{P}(f)=\E_\P(f)$ for the mathematical expectation of $f$ under $\P$.


As it is common in clustering,  we assume that the $k$-means live in a certain subset $C\subset \mathcal{B}$. We consider the collection of functions
\begin{equation}\label{Eqn:k-means-class}
    \mathcal{F}_{V_{k}(C)}=\left\{ f_a  :  a=\left(a_{1},\ldots,a_{k}\right)\in V_{k}(C)\right\},
\end{equation}
where $f_a:\mathcal{B}\to \R$ is defined by
\begin{equation}\label{Eqn:k-means-functions}
    f_a(z)=\underset{i=1,\ldots,k}{\operatorname{min}}\,\left(\left\|z-a_{i}\right\|_{\mathcal{B}}^{2}\right),\quad z\in \mathcal{B},
\end{equation}
and
\begin{equation}\label{Eqn:k-means-set}
V_{k}(C)=\left\{\left(a_{1},\ldots,a_{k}\right)\in C^{k}: a_{i}\neq a_{j} \text{ for }i\neq j\right\}.
\end{equation}


When the search for the $k$-means is restricted to $C\subset\mathcal{B}$, the problem \eqref{kmeans-pop} amounts to finding $f_\mu\in\mathcal{F}_{V_{k}(C)}$ (or, equivalently, the elements $\mu=(\mu_1,\dots,\mu_k)\in V_{k}(C)$) such that ${\operatorname{P}}(f_\mu)=W(k,C)$, where
\begin{equation}\label{WCSS-population-restricted}
W(k,C)= \inf_{f\in \mathcal{F}_{V_{k}(C)}} \P(f)
\end{equation}
is the population WCSS restricted to the set $C$. In other words, $\left\{\mu_{1},\ldots,\mu_{k}\right\}\in \mathcal{S}_{\P}(k)$ (an element of the set of all $k$-means) with $\mu_i\in C$ if and only if  the function $f_\mu$ in \eqref{Eqn:k-means-functions} with $\mu=\left(\mu_{1},\ldots,\mu_{k}\right)\in V_{k}(C)$ is a risk-minimizer of the operator ${\operatorname{P}}$ within the class $\mathcal{F}_{V_{k}(C)}$. Note that the parametrization of $\mathcal{F}_{V_{k}(C)}$ in terms of the set $V_{k}(C)$ in \eqref{Eqn:k-means-set} is not injective: different permutations of $a\in V_{k}(C)$ lead to the same function $f_a\in \mathcal{F}_{V_{k}(C)}$.

In the following theorem we calculate the asymptotic distribution of the empirical WCSS restricted to $C$, that is,
\begin{equation}\label{WCSS-empirical-restricted}
W_n(k,C) = \underset{f\in\mathcal{F}_{V_{k}(C)}}{\operatorname{inf}}\,\mathbb{P}_n(f).
\end{equation}
This result is key to characterize the uniqueness of the set of $k$-means and derive the test for $k$-means uniqueness in Section \ref{sec:test}. When dealing with risk minimization, it is often required that the class of functions in use is Glivenko-Cantelli (see, e.g.,  \cite[Section 2.1]{van der Vaart-Wellner-2023}) to ensure that ${\mathbb{P}_n}(f)$ and ${\operatorname{P}}(f)$ are close to each other uniformly over the class of functions. However, to derive the asymptotic distribution of $W_n(k,C)$ in \eqref{WCSS-population-restricted}  it is essential that $\mathcal{F}_{V_{k}(C)}$ satisfies the central limit theorem --i.e., $\mathcal{F}_{V_{k}(C)}$ has to be a \textit{Donsker class} (see the definition below)--, which is obviously a more demanding condition. Donsker classes are rather general and they have been already considered in this setting; see \cite{Caponnetto-Rakhlin-2006}.



Given a collection of real functions $\mathcal{F}$, $\ell^\infty(\mathcal{F})$ is the real Banach space of bounded functionals $\phi:\mathcal{F} \to \R$, equipped with the supremum norm, $\Vert \phi \Vert_\infty=\sup_{f\in \mathcal{F}}|\phi(f)|$. We recall that the class  $\mathcal{F}$ is \textit{$\P$-Donsker} if $\mathbb{G}_{n} \rightsquigarrow \mathbb{G}_\P$ in $\ell^\infty(\F)$, where `$\rightsquigarrow$' denotes weak convergence, $\mathbb{G}_{n}$ is the $\mathcal{F}$-indexed empirical process,
\begin{equation*}
\mathbb{G}_{n}(f)=n^{-1/2} \sum_{i=1}^n (f(X_i)-\P(f)),\quad f\in\mathcal{F},
\end{equation*}
and $\mathbb{G}_\P$ is a $\P$-Brownian bridge indexed by $\mathcal{F}$, that is, $\{ \mathbb{G}_\P (f) : f\in  \F \}$ is a zero-mean Gaussian process with covariance function
$\text{Cov}\left(\mathbb{G}_\P (f_1) , \mathbb{G}_\P (f_2)\right) =\P (f_1 f_2)-\P(f_1) \P (f_2)$, for $ f_1,f_2\in \F$. The class $\mathcal{F}$ is \textit{universal Donsker} if it is {$\P$-Donsker} for every probability measure $\P$ on the sample space.


The following result provides the asymptotic distribution of the normalized version of $W_n(k,C)$ under very general conditions.

\begin{theorem}\label{th:asymp}
Let $C\subset\mathcal{B}$ and let us consider  $W(k,C)$ and $W_n(k,C)$ the WCSS defined in \eqref{WCSS-population-restricted} and \eqref{WCSS-empirical-restricted}, respectively. Let us assume that $X$ fulfills \ref{itm:Int} and additionally that
the following two conditions are satisfied. 
\begin{description}
\labitem{(Bnd)}{itm:Bounded} \textit{Boundedness assumption.} The set $C$ is bounded in $\mathcal{B}$. In other words, we restrict the search for the $k$-means to a bounded set of the space.

\labitem{(Dnk)}{itm:Donsker} \textit{Donsker assumption.} The class $\mathcal{F}_{V_{k}(C)}$ in \eqref{Eqn:k-means-class} is $\P$-Donsker.
\end{description}

Then, we have that
\begin{equation}\label{Eqn:TnLimit}
T_{n}(k) =  \sqrt{n}\,\left(    W_n(k,C)  -  W(k,C) )      \right)
\cd T(k)=\underset{f\in  {S}_\P(k, C) }{\operatorname{inf}}\, {\mathbb{G}_{\operatorname{P}} (f)},
\end{equation}
where
\begin{equation}\label{completion-F}
{S}_\P(k, C)  = \left\{  f\in  \overline{\mathcal{F}}_{V_{k}(C)}  :  \P (f) =  W(k,C)  \right\}
\end{equation}
is the set of all minimizers of $\P$ over the class $\overline{\mathcal{F}}_{V_{k}(C)}\equiv$ the completion of $\mathcal{F}_{V_{k}(C)}$ with respect to the metric $d_{\mathrm{L}^2(\P)}$ with $d_{\mathrm{L}^2(\P)}^2(f,g)=\P(f-g)^2$ ($f,g\in \mathcal{F}_{V_{k}(C)}$).

				
\end{theorem}
\begin{proof}
First, by \ref{itm:Bounded}, $M=\sup_{b\in C} \Vert b \Vert_{\mathcal{B}}^2< \infty$. For any $f_a\in \mathcal{F}_{V_{k}(C)}$ with $a=(a_1,\dots,a_k)\in V_k(C)$, by \ref{itm:Int} we have that
\begin{equation}\label{bounded-L1}
\begin{split}
\operatorname{P}(f_a)&=\int_{\mathcal{B}}\underset{i=1,\ldots,k}{\operatorname{min}}\,\left(\left\|z-a_{i}\right\|_{\mathcal{B}}^{2}\right)\,\operatorname{dP}(z)\\
&=2^{2}\,\int_{\mathcal{B}}\underset{i=1,\ldots,k}{\operatorname{min}}\,\left(\left\|\frac{z-a_{i}}{2}\right\|_{\mathcal{B}}^{2}\right)\,\operatorname{dP}(z)\\
&\overset{(\ast)}{\leq}2^{2}\,\int_{\mathcal{B}}\underset{i=1,\ldots,k}{\operatorname{min}}\,\left(\frac{1}{2}\,\left(\left\|z\right\|_{\mathcal{B}}^{2}+\left\|a_{i}\right\|_{\mathcal{B}}^{2}\right)\right)\,\operatorname{dP}(z)\\
& \le  2  \,  (\E\Vert X\Vert_{\mathcal{B}}^2 + M) <\infty,
\end{split}
\end{equation}
where $(\ast)$ follows by the convexity of the norm and the square. Therefore, the functional $\P:\mathcal{F}_{V_{k}(C)} \to \R $ defined by $\P(f)=\E_\P(f)$ belongs to $\ell^\infty(\mathcal{F}_{V_{k}(C)})$.

Now, we use similar ideas as those in the proof of \cite[Theorem 6.1]{car20}. As $\mathcal{F}_{V_{k}(C)}$ is $\P$-Donsker, the space $(\mathcal{F}_{V_{k}(C)},d_\P)$ is totally bounded, where $d_\P$ is the pseudo-metric given by $d_{\text{P}}^2(f,g):= {\text{P} (f-g)^2-(\text{P}(f-g))^2}$, for $f,g\in \mathcal{F}_{V_{k}(C)}$; see \citet[Remark 3.7.27]{Gine-Nickl}. Also, the $\P$-Brownian bridge  $\mathbb{G}_\P\in \mathcal{C}_u(\mathcal{F}_{V_{k}(C)},d_\P)$ (uniformly continuous) a.s. By \eqref{bounded-L1}, the class $\mathcal{F}_{V_{k}(C)}$ is bounded in $\mathrm{L}^1(\P)$, i.e., $\sup_{f\in\mathcal{F}_{V_{k}(C)}} | \P (f)| <\infty$. This condition, joint to the fact that $(\mathcal{F}_{V_{k}(C)},d_\P)$ is totally bounded, implies that $(\mathcal{F}_{V_{k}(C)},d_{\mathrm{L}^2(\P)})$ is totally bounded, which follows from the same ideas as in the proof of \citet[Theorem 3.7.40, p. 262]{Gine-Nickl}. Note that the trajectories of $\mathbb{G}_\P$ also belong to $\mathcal{C}_u(\mathcal{F}_{V_{k}(C)},d_{\mathrm{L}^2(\P)})$ with probability 1 because $d_\P\le d_{\mathrm{L}^2(\P)}$. Finally, we can apply the extended delta method for the infimum in \cite[Corollary 2.5 (b)]{car20} to derive the asymptotic result in \eqref{Eqn:TnLimit}.
\end{proof}

We observe that the set of functions $\overline{\mathcal{F}}_{V_{k}(C)}$ in Theorem \ref{th:asymp} includes the set of all possible limits (in the completion of $\mathcal{F}_{V_{k}(C)}$ in \eqref{Eqn:k-means-class} with respect to $d_{\mathrm{L}^2(\P)}$) of minimizing sequences of the WCSS. Hence, ${S}_\P(k, C)$ in \eqref{completion-F} is the collection of sets of $k$-means (of the measure $\P$) restricted to $C$. Note that we need to complete $\mathcal{F}_{V_{k}(C)}$ to ensure the existence of minimizers. Finally, assumption \ref{itm:Donsker} is necessary to apply the uniform central limit theorem over $\mathcal{F}_{V_{k}(C)}$ and find the limit of the statistic $T_n(k)$ in \eqref{Eqn:TnLimit} by using the directional Hadamard differentiability of the infimum and the results in \cite{car20}.

\subsection{Donsker theorems for $k$-means}

In this section we establish conditions guaranteeing that the class of functions $\mathcal{F}_{V_{k}(C)}$ in \eqref{Eqn:k-means-class} is $\P$-Donsker, and hence the key assumption \ref{itm:Donsker} in Theorem \ref{th:asymp} is fulfilled.


In the classical, finite dimensional, setting of the $k$-means problem; when $X$ takes values in a ball of $\R^d$ (with the Euclidean distance), then $\mathcal{F}_{V_{k}(C)}$ in \eqref{Eqn:k-means-class} is universal Donsker; see \citet[Lemma 3.1]{Caponnetto-Rakhlin-2006}. This boundedness assumption is usually imposed in works related to $k$-means. In \cite{Telgarsky-Dasgupta-2013} some additional results are provided for distributions with $p\ge 4$ finite moments.


In a finite-dimensional context, the search of the $k$-means within a bounded set seems quite reasonable. However, in the infinite dimensional (functional) setting more demanding conditions are needed by the loss of the Heine-Borel's property. Sufficient conditions to ensure that a class of real functions $\mathcal{F}$ satisfies the Donsker property are usually related to the size of its coverings.


Given $C\subset\mathcal{B}$, it is said that $C$ is \textit{totally bounded} if for every $\varepsilon>0$ there exists $N\in\mathbb{N}$ and $x_{1},\ldots,x_{N}\in\mathcal{B}$ such that $C\subseteq\bigcup_{i=1}^{N}B_{d_\mathcal B}\left(x_{i},\varepsilon\right)$, where $B_{d_\mathcal B}$ stands for the ball with respect to the metric $d_\mathcal B$ induced by the norm in ${\mathcal B}$. Then, we say that $C$ is covered by $\left\{B_{d_\mathcal B}\left(x_{i},\varepsilon\right)\right\}_{i=1}^{N}$. For $\varepsilon>0$, let $N(\varepsilon,C,d_\mathcal B)$ be the \textit{$\varepsilon$-covering number} of $C$, that is, the minimal number of balls of radius $\varepsilon$ required to cover $C$. In the sequel  $N(\varepsilon,C,d)$ or $N(\varepsilon,C,\Vert\cdot\Vert)$ is the covering number associated with a general distance $d$ or a norm $\Vert\cdot\Vert$, respectively. 

Alternatively, if $\mathcal{F}$ is endowed with a metric $d$, it can be covered with brackets. For $l,u\in\mathcal{F}$, the \textit{bracket} is defined by
$$[l,u]=\{f\in\mathcal{F}:l\leq f\leq u\}.$$
An $\varepsilon$-bracket (with respect to $d$) is a bracket where $d(l,u)<\varepsilon$. The \textit{$\varepsilon$-bracketing number} $N_{[\,]}(\varepsilon,\mathcal{F},d)$ is the minimal number of $\varepsilon$-brackets necessary to cover $C$. The Donsker property is typically ensured for a class $\mathcal{F}$ through conditions that are associated with the covering and/or bracketing numbers; see \cite[Section 2.5]{van der Vaart-Wellner-2023}.


The following lemma links the bracketing numbers of $\mathcal{F}_{V_{k}(C)}$ in \eqref{Eqn:k-means-class} to the covering numbers of the underlying set $C$. We consider $\operatorname{L}^{2}(\operatorname{P})$, the usual $\operatorname{L}^{2}$-space with respect to the probability measure $\P$ induced by $X$.

\begin{lemma}\label{Lemma:Bracket-CoveringNumbers}
Let $C\subset\mathcal{B}$ be a totally bounded set. For the class $\mathcal{F}_{V_{k}(C)}$ in \eqref{Eqn:k-means-class}, we have that
\begin{equation}\label{Eqn:Lipschitz-bound}
		N_{[\,]}\left(\varepsilon,\mathcal{F}_{V_{k}(C)},\|\cdot\|_{\operatorname{L}^{2}(\operatorname{P})}\right)\leq N\left(\frac{\varepsilon}{\Delta},C,d_{\mathcal{B}}\right)^{k},
	\end{equation}
	where $\Delta=2\,\left(\int_{\mathcal{B}}\|z\|_{\mathcal{B}}^{2}\,\operatorname{dP}(z)+\underset{b\in C}{\operatorname{sup}}\,\left(\left\|b\right\|_{\mathcal{B}}\right)\right)$.
\end{lemma}
\begin{proof}
Let us first prove the following Lipschitz property for the functions in $\mathcal{F}_{V_{k}(C)}$:
\begin{equation}\label{Lipschitz-property}
\left|f_{a}(z)-f_{b}(z)\right|\leq2\,\left(\|z\|_{\mathcal{B}}+\underset{b\in C}{\operatorname{sup}}\,\left(\left\|b\right\|_{\mathcal{B}}\right)\right)\,d_{\infty}\left(a,b\right),
\end{equation}
 where $a=\left(a_{1},\ldots,a_{k}\right)$, $b=\left(b_{1},\ldots,b_{k}\right)\in V_{k}(C)$ and $d_{\infty}\left(a,b\right)=\underset{i=1,\ldots,k}{\operatorname{max}}\,\left(\left\|a_{i}-b_{i}\right\|_{\mathcal{B}}\right)$.
Indeed, we have that
\begin{equation*}
		\begin{split}
			\left|f_{a}(z)-f_{b}(z)\right|&=\left|\underset{i=1,\ldots,k}{\operatorname{min}}\,\left(\left\|z-a_{i}\right\|_{\mathcal{B}}^{2}\right)-\underset{i=1,\ldots,k}{\operatorname{min}}\,\left(\left\|z-b_{i}\right\|_{\mathcal{B}}^{2}\right)\right|\\
			&=\left|\underset{i=1,\ldots,k}{\operatorname{max}}\,\left(-\left\|z-a_{i}\right\|_{\mathcal{B}}^{2}\right)-\underset{i=1,\ldots,k}{\operatorname{max}}\,\left(-\left\|z-b_{i}\right\|_{\mathcal{B}}^{2}\right)\right|\\
			&\leq\underset{i=1,\ldots,k}{\operatorname{max}}\,\left(\left|\left\|z-a_{i}\right\|_{\mathcal{B}}^{2}-\left\|z-b_{i}\right\|_{\mathcal{B}}^{2}\right|\right)\\
			&=\underset{i=1,\ldots,k}{\operatorname{max}}\,\left(\left|\left\|z-a_{i}\right\|_{\mathcal{B}}+\left\|z-b_{i}\right\|_{\mathcal{B}}\right|\,\left|\left\|z-a_{i}\right\|_{\mathcal{B}}-\left\|z-b_{i}\right\|_{\mathcal{B}}\right|\right)\\
			&\leq\left(2\,\|z\|_{\mathcal{B}}+\underset{i=1,\ldots,k}{\operatorname{max}}\,\left(\left\|a_{i}\right\|_{\mathcal{B}}+\left\|b_{i}\right\|_{\mathcal{B}}\right)\right)\,\underset{i=1,\ldots,k}{\operatorname{max}}\,\left(\left\|a_{i}-b_{i}\right\|_{\mathcal{B}}\right)\\
    &\leq 2\,\left(\|z\|_{\mathcal{B}}+\underset{b\in C}{\operatorname{sup}}\,\left(\left\|b\right\|_{\mathcal{B}}\right)\right) \, d_{\infty}\left(a,b\right).
		\end{split}
	\end{equation*}
Therefore, \eqref{Lipschitz-property} holds. Finally, by \eqref{Lipschitz-property} and since $C$ is totally bounded, we can apply \cite[Theorem 2.7.17]{van der Vaart-Wellner-2023} to derive \eqref{Eqn:Lipschitz-bound}.
\end{proof}
Sufficient conditions for a class to be $\operatorname{P}$-Donsker in \cite[Theorem 19.5]{van der Vaart} are quite standard in empirical processes literature. The following theorem allows us to express this condition in terms of the, more geometrically motivated, covering numbers of $C$.

\begin{theorem}\label{Th:Donsker-kmeans}
Assume that
\begin{equation}\label{Donsker-condition}
\int_{0}^{1}\sqrt{\operatorname{log}\,\left(N\left(\varepsilon,C,d_{\mathcal{B}}\right)\right)}\,\operatorname{d}\varepsilon<\infty.
\end{equation}
Then, $\mathcal{F}_{V_{k}(C)}$ in \eqref{Eqn:k-means-class} is a $\P$-Donsker class.
\end{theorem}
\begin{proof}
	By Lemma \ref{Lemma:Bracket-CoveringNumbers} we have that
    \begin{equation}\label{donsker-inequality}
\int_{0}^{1}\sqrt{\operatorname{log}\left(N_{[\,]}\left(\varepsilon,\mathcal{F}_{V_{k}(C)},\|\cdot\|_{\operatorname{L}^{2}(\operatorname{P})}\right)\right)}\,\operatorname{d}\varepsilon\le \sqrt{k}\,\int_{0}^{1}\sqrt{\operatorname{log}\left(N\left(\frac{\varepsilon}{\Delta},C,d_{\mathcal{B}}\right)\right)}\,\operatorname{d}\varepsilon.
	\end{equation}
	Hence, if the right-hand side in \eqref{donsker-inequality} is finite, by \cite[Theorem 19.5]{van der Vaart} the class $\mathcal{F}_{V_{k}(C)}$ endowed with the $\operatorname{L}^{2}(\operatorname{P})$ is Donsker.
\end{proof}
The previous theorem extends the findings in \cite{pol82}, where the asymptotic Gaussianity of the empirical $k$-means is established. Let us now comment on the crucial condition \eqref{Donsker-condition} concerning the set $C$. This requirement might appear restrictive at first glance. However, it is satisfied by a broad range of examples. To begin with, it holds for every finite-dimensional bounded set. Moreover, every Vapnik-\u{C}ervonenkis (VC) class satisfies this condition, with its entropy numbers growing polynomially; for this and related results, we refer to \cite[Section 2.6]{van der Vaart-Wellner-2023}.

In the realm of functional data analysis (FDA), there are also several examples in which $\mathcal{F}_{V_{k}(C)}$ in \eqref{Eqn:k-means-class} is Donsker. For instance, unit balls in $\alpha$-Hölder continuous, Sobolev, and Besov function spaces satisfy this condition under mild restrictions on the parameters; see \cite[Section 2.7]{van der Vaart-Wellner-2023}. It is worth noting that these sets are dense in the unit balls of the space of continuous and bounded functions 
and in $\operatorname{L}^{p}$ spaces, respectively. In essence, the preceding lines imply that there is no loss of generality by searching for sets of $k$-means within dense subsets of \enquote{smooth functions}, which is in fact a very common working practice in this area. Therefore, the results of this work can also be applied in FDA.

Other subsets of commonly used spaces also satisfy \ref{itm:Donsker}. The following corollary lists some of these examples.
\begin{corollary}
Assume one of the following conditions is satisfied:
\begin{enumerate}[(a)]
\item $C$ is a subset of a finite dimensional Banach space $\mathcal{B}$.

\item $C$ is a Vapnik-\u{C}ervonenkis (VC) class of functions endowed with the $\operatorname{L}^{p}(\operatorname{P})$-metric, with $p\in[1,\infty)$.

\item $C$ is the unit ball of the space $\alpha$-Hölder continuous functions $\mathcal{C}^{\alpha}(\mathcal{X})$ with $\alpha>0$, where $\mathcal{X}\subset\mathbb{R}^{l}$ bounded, convex and with no empty interior.

\item $C$ is the unit ball of the Sobolev space $W^{\alpha,p}(\mathcal{X})$, where $\mathcal{X}$ is a bounded Lipschitz domain of $\mathbb{R}^{l}$, endowed with $\operatorname{L}^{r}$-norm (with $r,p\in[1,\infty]$) and $\alpha>\operatorname{max}\left(0, l\left(\frac{1}{p}-\frac{1}{r}\right)\right)$. More generally, $C$ is a subset of the unit ball of the Besov space $\mathcal{B}_{p,q}^{\alpha}(\mathcal{X})$ (with $q\in[1,\infty]$) endowed with the $\operatorname{L}^{r}$-norm.

\item $C$ is the class of monotone functions on $[0,1]$ endowed with an $\operatorname{L}^{p}$-norm, with $p\in[1,\infty]$.

\item $C$ is the class of all convex functions on a compact convex subset $D\subset\mathbb{R}^{l}$ with values on $[0,1]$ endowed with the $\operatorname{L}^{p}$-norm with $p\in[1,\infty]$.
\end{enumerate}
Then $\mathcal{F}_{V_{k}(C)}$ in \eqref{Eqn:k-means-class} is a $\operatorname{P}$-Donsker class.
\end{corollary}
\begin{proof}
It is enough to check condition \eqref{Donsker-condition} and apply Theorem \ref{Th:Donsker-kmeans} for each of the considered scenarios. To see (a), assume without loss of generality that $C=\mathbb{R}^{l}$. Now, we identify every point $\alpha=(\alpha_1,\dots,\alpha_l)\in\mathbb{R}^{l}$ with the indicator function ${1}_{A_\alpha}$, where $A_\alpha=\left(-\infty,\alpha_{1}\right]\times\cdots\times\left(-\infty,\alpha_{l}\right]$. By \cite[Example 2.6.1]{van der Vaart-Wellner-2023}, $\mathcal{G}=\left\{ {1}_{A_\alpha} : \alpha\in\mathbb{R}^{l}\right\}$ is a VC-class of functions that can be endowed with $\operatorname{L}^{p}(\operatorname{P})$-metric, for $p\in[1,\infty)$. Hence (a) is a consequence of (b).

In part (b), the decreasing of the covering numbers is polynomial in $\frac{1}{\varepsilon}$; see
\cite[Theorem 2.6.7]{van der Vaart-Wellner-2023}.

In part (c), to check that the decreasing of the logarithm of covering numbers is polynomial in $\frac{1}{\varepsilon}$ we refer to \cite[Theorem 2.7.1]{van der Vaart-Wellner-2023}.

Part (d) follows from \cite[Theorem 2.7.4]{van der Vaart-Wellner-2023}, (e) from \cite[Theorem 2.7.9]{van der Vaart-Wellner-2023} and (f) from \cite[Theorem 2.7.14]{van der Vaart-Wellner-2023}.
\end{proof}


\section{A characterization and test of $k$-means uniqueness}\label{sec:test}

We establish here necessary and sufficient conditions for the uniqueness of the set of $k$-means within a subset $C\subset \mathcal{B}$. Subsequently, we use these ideas to construct a bootstrap test of $k$-means uniqueness.

For $C\subset \mathcal{B}$, we say that $\P$ satisfies the \textit{$k$-means uniqueness property in $C$}, $\text{UP}(k,C)$, or $\P \in \text{UP}(k,C)$, if the set ${S}_\P(k,C)$ in \eqref{completion-F} has cardinality one. The following result, a consequence of Theorem \ref{th:asymp}, characterizes this situation in terms of the properties of the limit variable $T(k)$ in \eqref{Eqn:TnLimit}.

\begin{corollary}\label{cor:thecorollary}
Under assumptions  \ref{itm:Int}, \ref{itm:Bounded} and \ref{itm:Donsker}, the following three assertions are equivalent.

\begin{enumerate}
\item[(i)] $\P \in \textrm{\rm UP}(k,C)$ (uniqueness of $k$-means on $C$).

\item[(ii)]  The variable $T(k)$ in \eqref{Eqn:TnLimit} is normally distributed with mean zero.

\item[(iii)] The variable $T(k)$ in \eqref{Eqn:TnLimit} has zero mean.

\end{enumerate}
\end{corollary}
\begin{proof}
If (i) holds, then there exists a unique minimizer $f^-\in \overline{\mathcal{F}}_{V_{k}(C)}$ such that ${S}_\P(k,C)=\{f^-\}$. From Theorem \ref{th:asymp}, we obtain that $T(k)={\mathbb{G}_{\operatorname{P}} (f^-)}$, which has normal distribution with mean zero. Therefore, (ii) is satisfied. The implication (ii) $\Rightarrow$ (iii) is direct. Finally, assume that (iii) holds. Observe that $T(k)\le_{\text{st}} {\mathbb{G}_{\operatorname{P}} (f)}$, for each $f\in {S}_\P(k, C)$, where `$\le_{\text{st}}$' stands for the usual stochastic order. Since $\E T(k) = 0 =\E \mathbb{G}_{\operatorname{P}} (f)$, we conclude that $T(k)=_{\text{st}} \mathbb{G}_{\operatorname{P}} (f)$; see \cite[Theorem 1.A.8.]{Shaked-2007}. Hence, we obtain that $T(k)$ is normally distributed. Finally, as
$T(k)=\underset{f\in  {S}_\P(k, C) }{\operatorname{inf}}\, {\mathbb{G}_{\operatorname{P}} (f)}$ is the infimum of normal variables, its distribution only can be normal when the infimum is taken over a set of cardinality one; if $f_1,f_2\in  {S}_\P(k, C) $ with $f_1\ne f_2$, then $\min\{\mathbb{G}_{\operatorname{P}}(f_1),\mathbb{G}_{\operatorname{P}}(f_2) \} $ is not normally distributed. Hence, we conclude that the cardinality of ${S}_\P(k, C)$ is necessarily one and (i) holds. 	
\end{proof}

Corollary \ref{cor:thecorollary} makes it possible to develop various discrepancy measures to perform a statistical test for the following hypothesis testing problem on the uniqueness of the set of $k$-means:
\begin{equation}\label{Test.UP(k)}
\begin{cases}
H_0: & \P \in \text{UP}(k,C)\qquad \text{(uniqueness of $k$-means in $C$)},\\
H_1: &  \P \notin \text{UP}(k,C)\qquad \text{(non-uniqueness of $k$-means in $C$).}
\end{cases}
\end{equation}

We note that $H_0$ is equivalent to the asymptotic distribution of $T_n(k)$ in \eqref{Eqn:TnLimit} being normal with zero mean. Therefore, we can use any normality test available in the literature if we have a sufficient number of data. Further, under $H_0$ we have that the expectation of $T(k)$ is zero, and $\E T(k)<0$ under $H_1$. In particular, the testing problem  in \eqref{Test.UP(k)} can be equivalently stated in either of the following two ways:
\begin{equation}\label{Test.UP(k)-2}
\text{(a)}\quad
\begin{cases}
H_0: & T(k)\sim\mathcal{N}(0, \sigma^2),\\
H_1: &  \text{no }H_0,
\end{cases}
\qquad\quad \text{(b)}\quad
\begin{cases}
H_0: & \E T(k)=0,\\
H_1: & \E T(k)<0,
\end{cases}
\end{equation}
where $T(k)$ is the limiting distribution of $T_n(k)$ in equation \eqref{Eqn:TnLimit}.

In practice, we usually observe one sample $x_{1},\ldots,x_{n}$ of $X$ and can only compute a single value \( W_n(k, C) \) of the WCSS in \eqref{WCSS-empirical-restricted}. Therefore, to carry out the uniqueness test \eqref{Test.UP(k)} we have to resort to a bootstrap scheme to approximate the limit distribution $T(k)$.  Specifically, to construct a test statistic for \eqref{Test.UP(k)-2} we propose the following bootstrap procedure.

\noindent \textbf{Step 1.} Extract \(B\) bootstrap samples from the original sample.
\begin{equation} \label{flowchartKMUtest}
    \begin{array}{cccc}
\mbox{Original sample} & & \mbox{Bootstrap sample} \\
\begin{array}{l}
x_{1},\ldots,x_{n}
\end{array}
&
\begin{array}{c}
\longrightarrow
\end{array}
&
\begin{array}{l}
x_{1}^{*b},\ldots,x_{n}^{*b},
\end{array}
&
b=1,\ldots,B.
\end{array}
\end{equation}

\noindent  \textbf{Step 2.} For the original sample and each bootstrap sample in \eqref{flowchartKMUtest}, compute the empirical WCSS and bootstrapped empirical WCSS, \( W_n(k, C) \) and \( W_n^{*b}(k, C) \) for \(b=1,\ldots,B\), respectively.

\noindent  \textbf{Step 3.} For each \(b=1,\ldots,B\), compute the bootstrapped version of $T_n(k)$:
\begin{equation}\label{bootstrapped statistic}
T_{n}^{*b}(k) =  \sqrt{n}\,\left(    W_n^{*b}(k,C)  -  W_n(k,C)\right).
\end{equation}

\noindent  \textbf{Step 4.} Use the values $T_{n}^{*b}(k)$ (\(b=1,\ldots,B\)) to carry out one of the tests in \eqref{Test.UP(k)-2}.

The consistency of this proposed bootstrap procedure rests on the approximation of the limiting distribution $T(k)$ by the values in \eqref{bootstrapped statistic}.
Hence, it is essential that the bootstrap estimator of this quantity is consistent. This result is a consequence of several previous works and we collect it in the following proposition where the a.s. consistency is understood as in \cite{van der Vaart-Wellner-2023}[Theorem 3.7.2  (iii)].

\begin{proposition}
Let us assume that \ref{itm:Bounded} and \ref{itm:Donsker} hold and $\E\Vert X\Vert_{\mathcal{B}}^4<\infty$. Under $H_0$ in \eqref{Test.UP(k)} (uniqueness of $k$-means), the bootstrap estimator of $T(k)$ in \eqref{bootstrapped statistic} is a.s. consistent.
\end{proposition}

\begin{proof}
We first note that $T_{n}(k)= \sqrt{n}\,\left(   \iota(\mathbb{P}_n)  -  \iota(\P)   \right)$, where the functional $\iota$ is given by $\iota(g)=\inf_{f\in \mathcal{F}_{V_{k}(C)}} g(f)$, for $g\in\ell^\infty(\mathcal{F}_{V_{k}(C)})$. To conclude that the bootstrap estimator of $T_{n}(k)$ in \eqref{bootstrapped statistic} is a.s. consistent is enough to check that we are in the conditions to apply  \cite[Theorem 3.1]{Fang-Santos}. To see this, we need to verify the following three points:
\begin{itemize}
    \item [(i)] The empirical bootstrap process is a.s. consistent. This result follows from \ref{itm:Donsker} and \cite{van der Vaart-Wellner-2023}[Theorem 3.7.2]. Note that the finite forth moment of $X$ together with triangle inequality and \ref{itm:Bounded} imply that $\P\Vert f-\P(f) \Vert_{\mathcal{F}_{V_{k}(C)}}^2<\infty$.

\item [(ii)] The limit of the underlying (empirical) process, $\sqrt{n}\,\left(   \mathbb{P}_n -  \P  \right)$ is Gaussian. This holds since $\mathcal{F}_{V_{k}(C)}$ in \eqref{Eqn:k-means-class} is Donsker and the limiting process is $\mathbb{G}_\P$, the Brownian bridge indexed by $\mathcal{F}_{V_{k}(C)}$.

\item [(ii)] The functional $\iota$ is fully Hadamard differentiable at $\P$. To see this, we note that under $H_0$ there exists a unique minimizer $f^-\in\overline{\mathcal{F}}_{V_{k}(C)}$ of the WCSS in $C$, that is, ${S}_\P(k, C) =\{ f^- \} $. Hence, by \cite[Corollary 2.5]{car20}, the Hadamard directional derivative of the infimum $\iota(g)=\inf_{f\in \mathcal{F}_{V_{k}(C)}} g(f)$ evaluated at $\P$ is $\iota'_{\P}(g)=g(f^-)$, for $g\in\ell^\infty(\mathcal{F}_{V_{k}(C)})$. This derivative is linear and hence, by \cite[Proposition 2.1]{Fang-Santos}, $\iota$ is fully Hadamard differentiable at $\P$.
\end{itemize}
Thus, the proof is complete.

\end{proof}

\section{Empirical results}\label{sec:emp}

The aim of this section is to provide some insights about the theoretical results established in this work. In particular, we illustrate by numerical simulations the practical performance of the  $k$-means uniqueness testing problem \eqref{Test.UP(k)} (by means of \eqref{Test.UP(k)-2}) proposed in Section \ref{sec:test}. We analyze in depth all the examples in Section \ref{sec:non-unique}. We thus consider the different scenarios discussed in this paper: $\text{UP}(k)$, uniqueness of the $k$-means set; $\text{CNU}(k)$ continuous non-uniqueness; and $\text{DNU}(k)$ discrete non-uniqueness. The $\text{CNU}(k)$ case could be seen as ``more pathological'' since it entails the existence of different sets of $k$-means arbitrarily close (in Hausdorff distance), which usually leads to algorithm instability problems.


The numerical experiments include:
\begin{itemize}
    \item[(a)]  Several simulations designed to evaluate the performance of the testing problem \eqref{Test.UP(k)} for $k$-means uniqueness. Specifically, we check the empirical level and power of the bootstrap test based on the values in \eqref{bootstrapped statistic}.

 \item[(b)] A record of the computation times of each experiment to identify differences between the analyzed models.

    \item[(c)] A toy example showing how the uniqueness of the set of $k$-means can be determined or ruled out if we know the distribution of the population by means of Monte Carlo simulations.

    \item[(d)] Graphical illustrations of the consistency of the empirical $k$-means in Theorem \ref{th_cluster_consistency} for the different cases of non-uniqueness.
\end{itemize}

\subsection{Models and technical details}

We consider all the models described in Section \ref{sec:non-unique} (Examples \ref{Example-uniform}-\ref{Example.last}). We recall that C$i$k$j$ stands for a model in which the data are drawn for a mixture of $i$ distributions and we decide to look for $k=j$ centers. Ideally, if we had prior information on the underlying distribution (which seldom happens in practice), the $k$-means parameter should be typically chosen to match the  number of mixture components, that is, $i=j$. As commented in Section \ref{sec:non-unique}, non-uniqueness is a possible side effect of a ``wrong'' choice of $k$ when compared with the true number of components of the mixture. In those models with Gaussian variables (Examples \ref{Example-C1k2}-\ref{Example.last}) we also simulate from distributions in dimensions $d=2, 6, 12, 18, 24$. This extension to higher dimensions is done by filling with zeros the mean vectors of each component and by taking identity covariance matrices of the corresponding dimension.

We have used the Hartigan-Wong algorithm as implemented in the \texttt{kmeans} command of \texttt{R}, selecting 20 random starts (by choosing \texttt{nstart=20}), thus taking into account the fact that the algorithm provides often just a local minimum.

In all experiments related to the uniqueness testing problem \eqref{Test.UP(k)} we have set the significance level at $\alpha=0.05$, and the number of bootstrap samples is $B=1000$. Taking a larger value of $B$ does not significantly change the results. We have also replicated each experiment 200 times to record the percentage of rejections of the null hypothesis (uniqueness) among the 200 replicates.

For the experiments, we have used a a 512-core computer Dual AMD MILAN 7453 Series Processors.
Though the computational cost of our test, for a given data set, is perfectly affordable by a standard computer, a simulation study as this one is much more demanding; just to give a hint, the whole study presented here (including the unavoidable associated debugging process) took more than eight months of intensive computation.

\subsection{Testing for uniqueness}\label{Subsec:Testing}


Before presenting the outputs for the hypothesis testing problem \eqref{Test.UP(k)}, let us describe first the procedure in detail. According to Corollary \ref{cor:thecorollary}, \eqref{Test.UP(k)} is equivalent to test that $T_n(k)$ in \eqref{Eqn:TnLimit} is asymptotically normal which is also equivalent to the fact that the expected value of $T(k)$ in \eqref{Eqn:TnLimit} is zero. We have found that testing problem for the mean in \eqref {Test.UP(k)-2} (b) consistently provides more powerful results than the different normality tests that we have tried (Kolmogorov, Anderson-Darling, Jarque-Bera, etc.). Therefore, we have used \eqref {Test.UP(k)-2} (b) to solve the problem of interest in \eqref{Test.UP(k)}.

\

Specifically, given an observed sample $x_{1},\ldots,x_{n}$ from $\operatorname{P}$, we proceed as follows:

\noindent \textbf{Step 1.} As described in Section~\ref{sec:test}, compute the values $T_{n}^{*1}(k),\dots, T_{n}^{*B}(k)$ (with $B=1000$), where $T_{n}^{*b}(k)$ is in \eqref{bootstrapped statistic}.

\noindent \textbf{Step 2.} Calculate the standardized mean of $T_{n}^{*1}(k),\dots, T_{n}^{*B}(k)$, that is,
$$ \bar{T}^*(k)  = \frac{1}{s_n^* \sqrt{B}} \sum_{b=1}^B T_{n}^{*b}(k), $$
where $s_n^*$ is the (quasi)standard deviation of $T_{n}^{*1}(k),\dots, T_{n}^{*B}(k)$.

\noindent \textbf{Step 3.} We reject the null hypothesis in \eqref{Test.UP(k)} if $\bar{T}^*(k)<q(\alpha)$, the $\alpha$-quantile of a standard normal distribution.


The models are described in Section \ref{sec:non-unique}, Examples \ref{Example-C1k2}-\ref{Example.last}. Recall that models C2k2-1, C2k2-2, C2k2-3, C3k3 and C3k2 are under $H_0$ (uniqueness) and models C1k2, TC3k2, C2k3 are under $H_1$ (non-uniqueness). Further, models C1k2 and C2k3 corresponds to the case $\text{CNU}(2)$ and $\text{CNU}(3)$, respectively, and TC3k2 to $\text{DNU}(2)$.

We have designed the experiments to see the effect of the sample size $n$ on the test (control of the significance level and power). We have considered the values
$$n=1000, 2000, 5000, 10000, 20000, 50000, 100000, 200000, 500000.$$
We also want to check how the dimension of the observed vectors affects the outputs by considering $d=2,6,12,18,24$. The processing time of each experiment has been saved to study the computational cost of the algorithm for the different cases of non-uniqueness. Due to the large number of outputs obtained, we will summarize the main findings in a few graphs. However, other complementary results can be found in the Supplementary Material file.

\begin{figure}[h]
\centering
    \includegraphics[width=1 \textwidth]{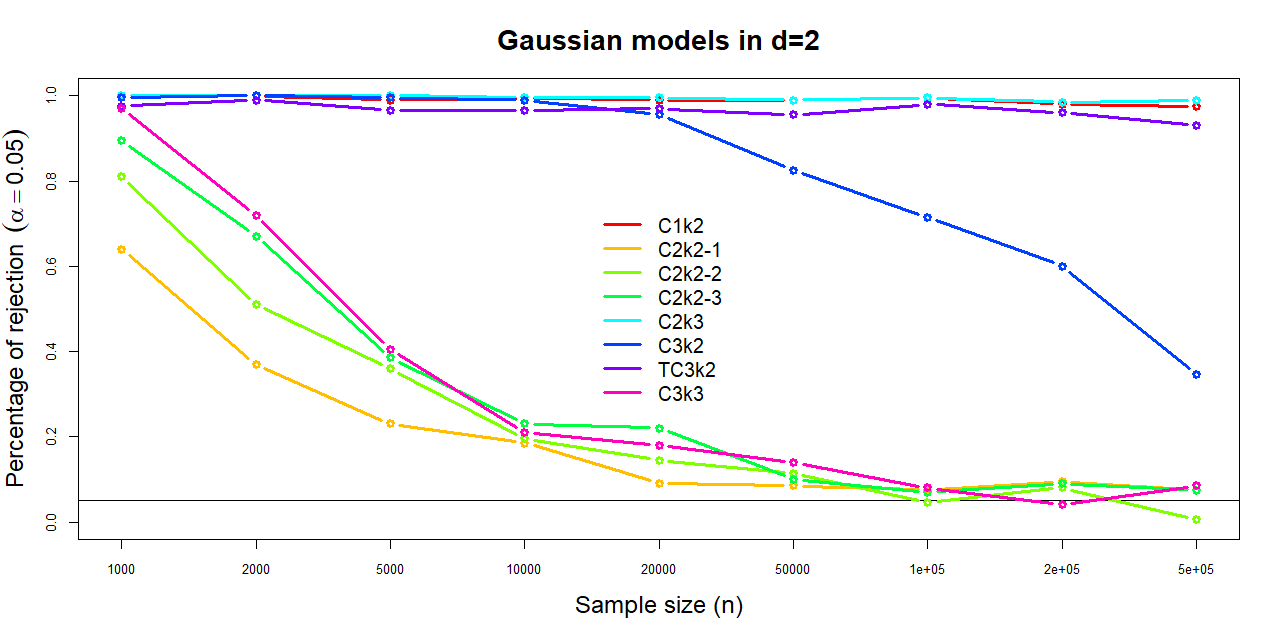}
    \caption{Percentage of rejections of the hypothesis of uniqueness in 8 models as a function of the sample size in dimension $d=2$. The thin black horizontal line corresponds to significance level $0.05$.
    }\label{fig:d2}
\end{figure}

\begin{figure}[h]\label{fig:d24}
\centering
    \includegraphics[width=1 \textwidth]{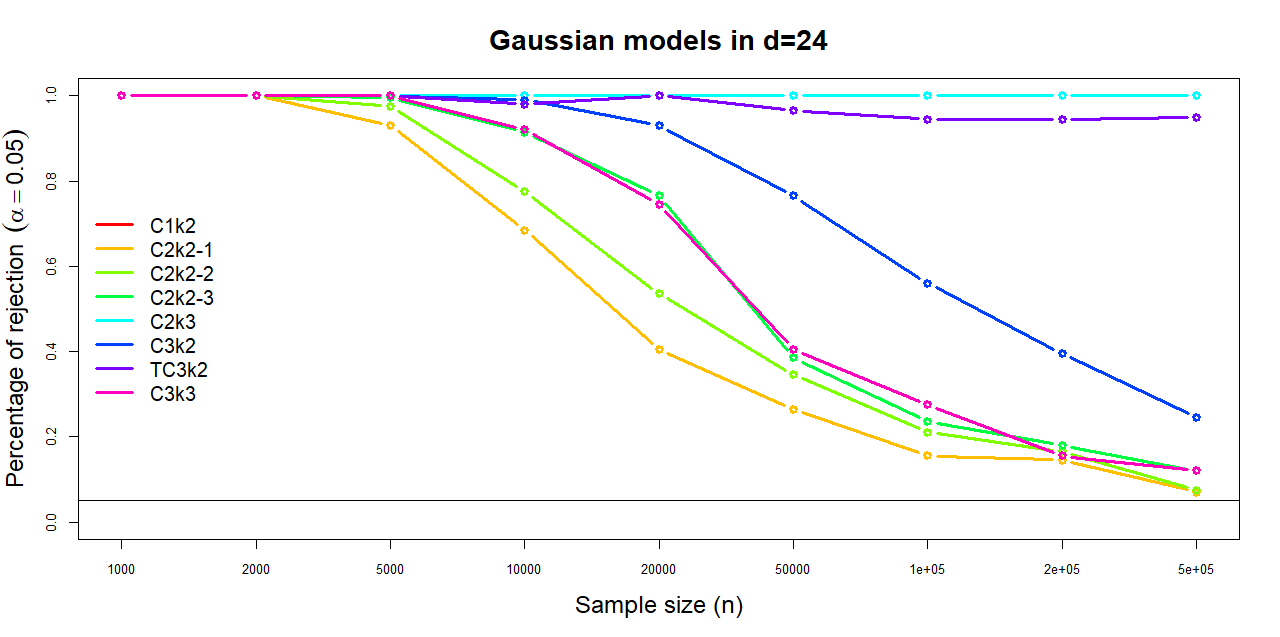}
    \caption{Percentage of rejections of the hypothesis of uniqueness in 8 models as a function of the sample size in dimension $d=24$. The thin black horizontal line corresponds to significance level $0.05$.
    }\label{fig:d24}
\end{figure}

In Figures \ref{fig:d2} and \ref{fig:d24} we see the percentage of rejections (among the 200 repetitions of each experiment) of the 8 considered models, with a fixed $\alpha=0.05$, of our proposal in dimension $d=2$ and $d=24$, respectively. The intermediate dimensions behave similarly and the corresponding graphs are included in the Supplementary Material.

We first note that Type I error is not controlled up to very large sample sizes. The test tends to over-reject under the null hypothesis systematically with lower sample sizes. This might be due in part to the fact that the procedure is based on the limiting behavior of the WCSS and that the  problem itself is very complex, requiring a lot of sample information for acceptable performance. As the sample size $n$ increases, the percentage of rejection of $H_0$ approaches the nominal value $\alpha=0.05$ in the uniqueness models C2k2-1, C2k2-2, C2k2-3, and C3k3. It is remarkable that the test detects the C2k2-3 case since it is a mixture of two very similar Gaussian distributions; see Example \ref{Example-C2k2} and Figure \ref{fig: C2k2-Mofified} (c) and (d). The C3k2 model takes the longest time to lower the rejection rate maybe because the choice of $k$ is wrong and this leads to a higher instability of the algorithm. In fact, a similar model with less variability would behave like Example \ref{Example-uniform} with small $r$ which is a case of non-uniqueness. Thus, we see that the test can detect cases of a bad choice of $k$ as the model C3k2 since it rejects $H_0$ in this case until very large sample sizes. All models under $H_1$ (continuous and discrete non-uniqueness) are successfully detected. As the dimension increases, even larger sample sizes are needed due to the noise of the added variables and perhaps because at higher dimensions there are more local optima in which the algorithm can get stuck.


\begin{figure}[h]
\centering
    \includegraphics[width=1 \textwidth]{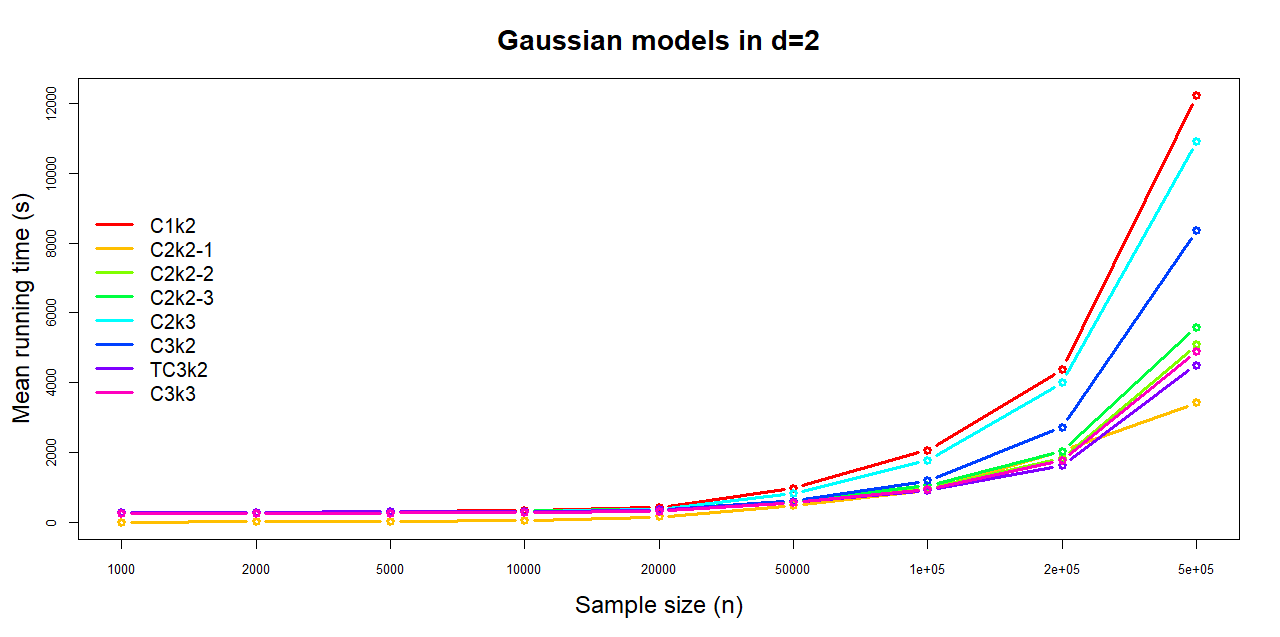}
    \caption{Average execution times (in seconds) of the $k$-means algorithm for the 200 replicates of each experiment in dimension $d=2$.
    }\label{fig:d2-time}
\end{figure}

\begin{figure}[h]
\centering
    \includegraphics[width=1 \textwidth]{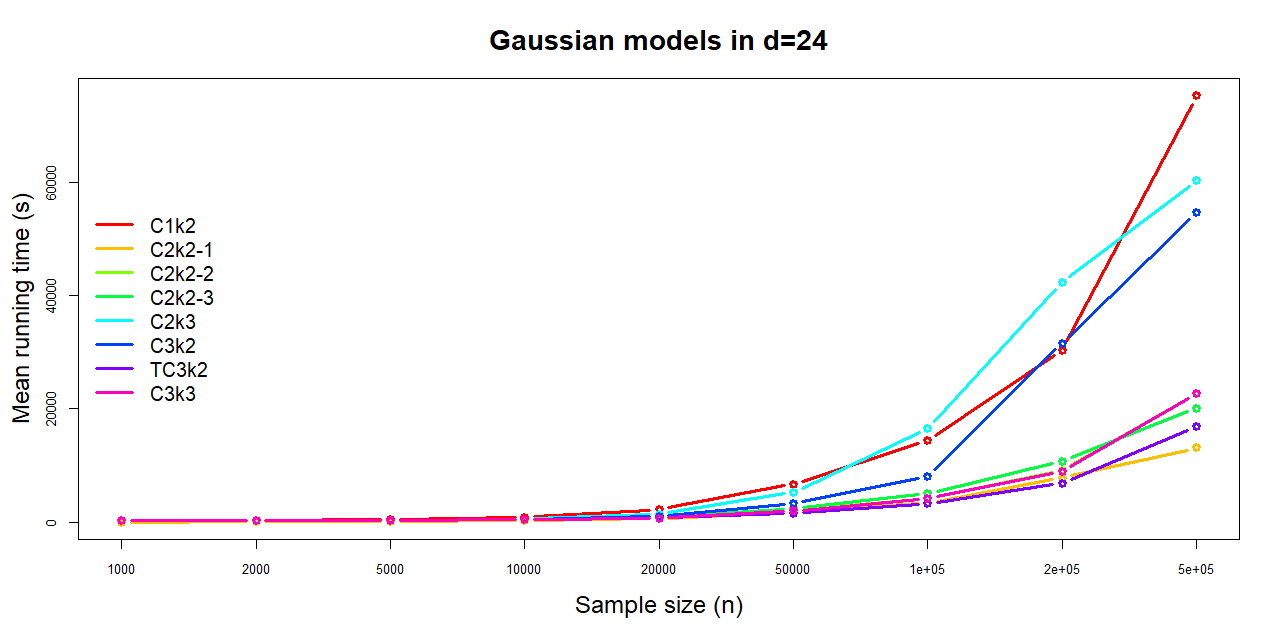}
    \caption{Average execution times (in seconds) of the k-means algorithm for the 200 replicates of each experiment in dimension $d=24$.
    }\label{fig:d24-time}
\end{figure}

Figures \ref{fig:d2-time} and \ref{fig:d24-time} show the mean execution times of the $k$-means algorithm for the 200 replicates of each experiment in dimension $d=2$ and $d=24$, respectively. We observe that the processing times of the continuous non-uniqueness cases (models C1k2 and C2k3) are significantly higher than the rest, especially for high sample sizes and high dimension. The next longest running model is the case C3k2, which corresponds to a bad choice of $k$ and is very close to non-uniqueness. It is important to note that the TC3k2 model (discrete non-uniqueness) has very similar execution times to the models in which there is uniqueness with a good choice of $k$. This supports our assumption that the instability of the algorithms occurs mainly in the $\text{CNU}(k)$ case. Information about running times in intermediate dimensions is collected in the Supplementary Material.

\subsection{Assessing uniqueness by Monte Carlo simulations}

As commented above, even if we know the true distribution of the population, it might be quite difficult  to determine whether the uniqueness property $\text{UP}(k)$ is satisfied or not. Nevertheless, as a further application of Corollary \ref{cor:thecorollary}, we could check this property using Monte Carlo simulations. Given a probability measure $\P$, we first note that in general it is not possible to use a standard Monte Carlo approach in this problem. This is because it is not easy to calculate analytically the population WCSS appearing in \eqref{Eqn:TnLimit} to approximate the distribution of $T(k)$ in Corollary \ref{cor:thecorollary}. Thus, in practice we can draw a large sample from $\P$ and use the same bootstrap scheme described in the previous Section \ref{Subsec:Testing} to perform the uniqueness test \eqref{Test.UP(k)}. If we repeat this procedure a sufficient number of times we will be able to identify whether $\P\in \text{UP}(k)$ or not.

\begin{figure}[h]
\centering
    \includegraphics[width=1 \textwidth]{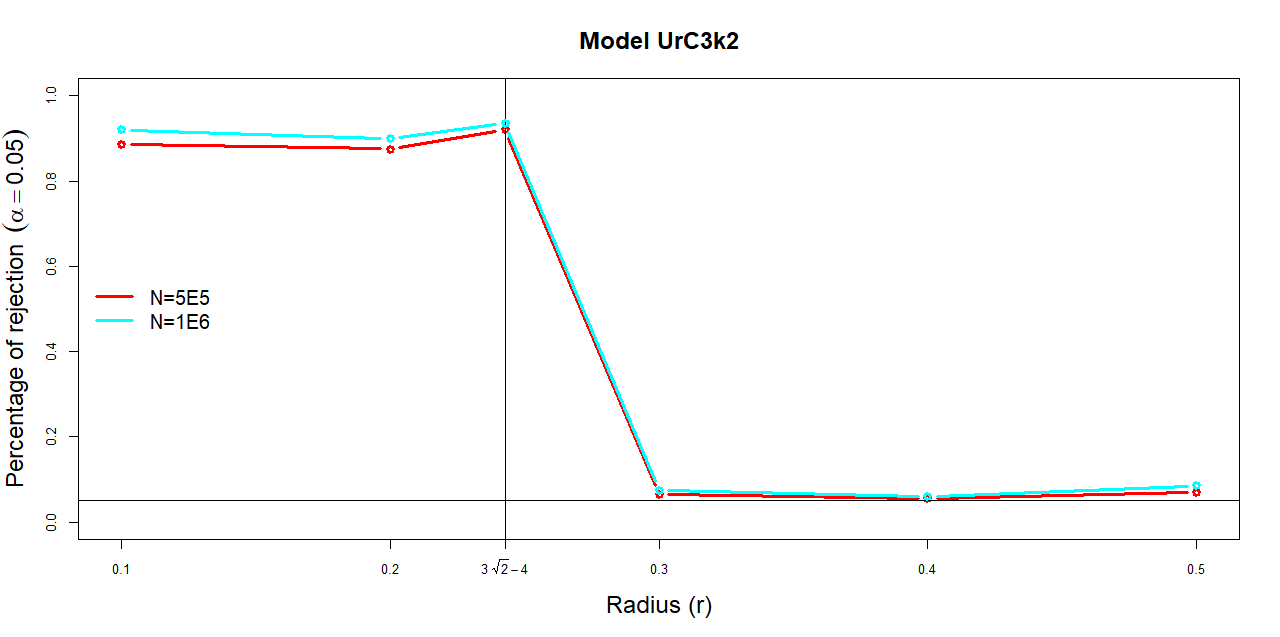}
    \caption{Percentage of rejections (among 200 replicates) in Example \ref{Example-uniform} of the hypothesis of uniqueness based of the procedure described in Section \ref{Subsec:Testing} for samples of sizes $n=5\times 10^5$ and $n=  10^6$.}
    \label{fig:Testing-U}
\end{figure}

To illustrate this idea we consider the distributions in Example \ref{Example-uniform}. In this case we have determined when $k$-means uniqueness holds (for $k=2$), as well as the population WCSS, as a function of the value of the radius $r$ of the underlying uniform variables. In Figure \ref{fig:Testing-U} we see the percentage of rejections (among 200 replicates) of the bootstrap test for uniqueness proposed in Section \ref{Subsec:Testing}, where we simulate $n=5\times 10^5$ and $n=10^6$ data from $\P_r$ in Example \ref{Example-uniform} with radii $r=0.1, 0.2, 3\sqrt{2}-4, 0.3, 0.4, 0.5$. We observe how the procedure perfectly detects cases of non-uniqueness ($r\le 3\sqrt{2}-4$) and uniqueness ($r>3\sqrt{2}-4$).

\begin{figure}[h]
\centering
    \includegraphics[width=1 \textwidth]{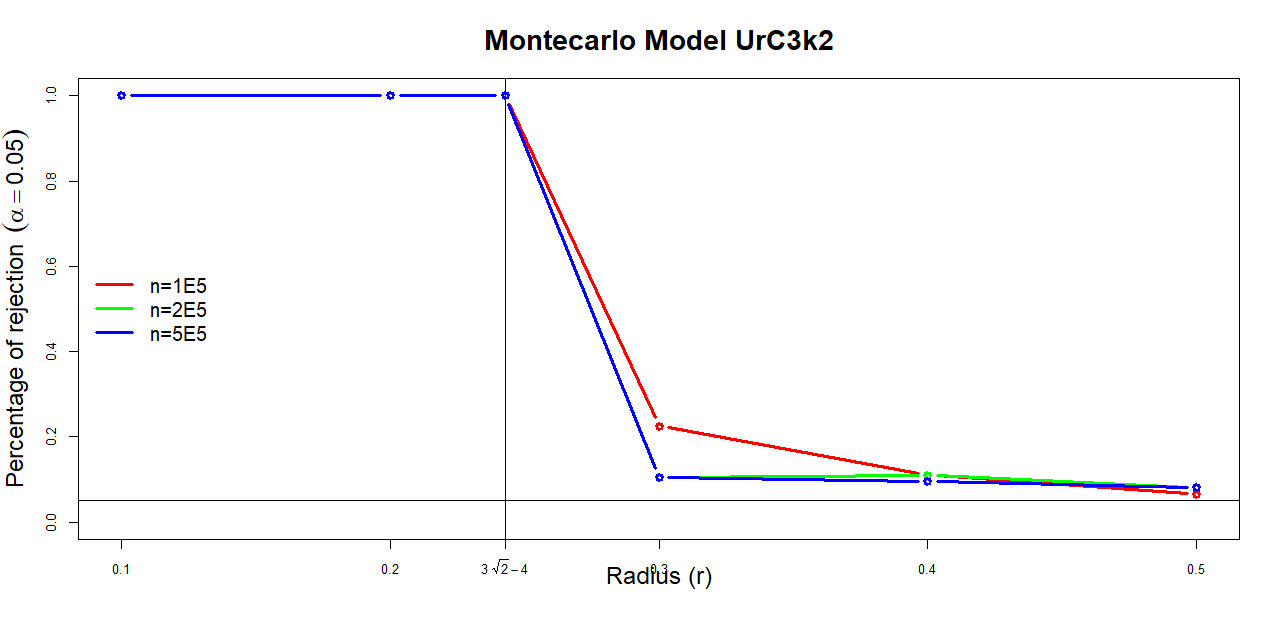}
    \caption{Percentage of rejections (among 200 replicates) in Example \ref{Example-uniform} of the hypothesis of uniqueness using the population WCSS for samples of sizes $n=10^5, 2\times 10^5, 5\times 10^5$.}\label{fig:Testing-MC}
\end{figure}

In Figure \ref{fig:Testing-MC} we have used the population WCSS computed in Example \ref{Example-uniform} to carry out a true Monte Carlo simulation by taking 200 samples of sizes $n=10^5, 2\times 10^5, 5\times 10^5$. It is interesting to note that the results of the MC simulations (using the true population WCSS value) are very similar to those of the bootstrap test.

\subsection{Some graphical illustrations of Theorem \ref{th_cluster_consistency}}

\begin{figure}[h]
\centering
    \begin{tabular}{cc}
     \includegraphics[width=0.38\textwidth]{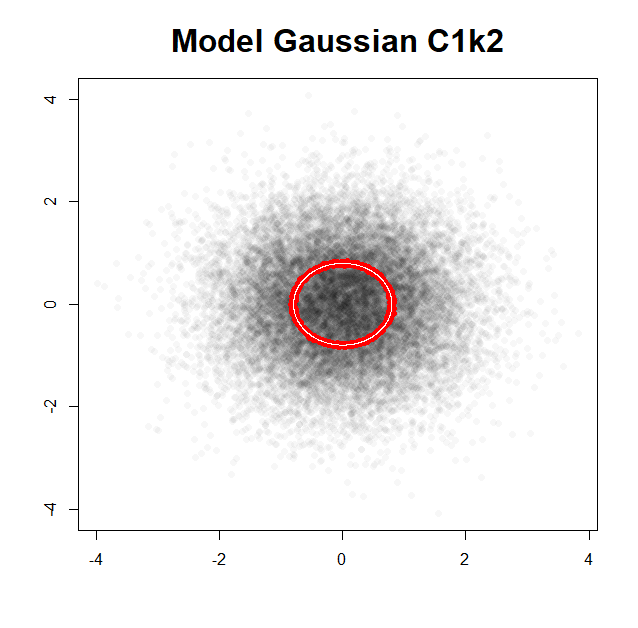} &
    \includegraphics[width=0.62\textwidth]{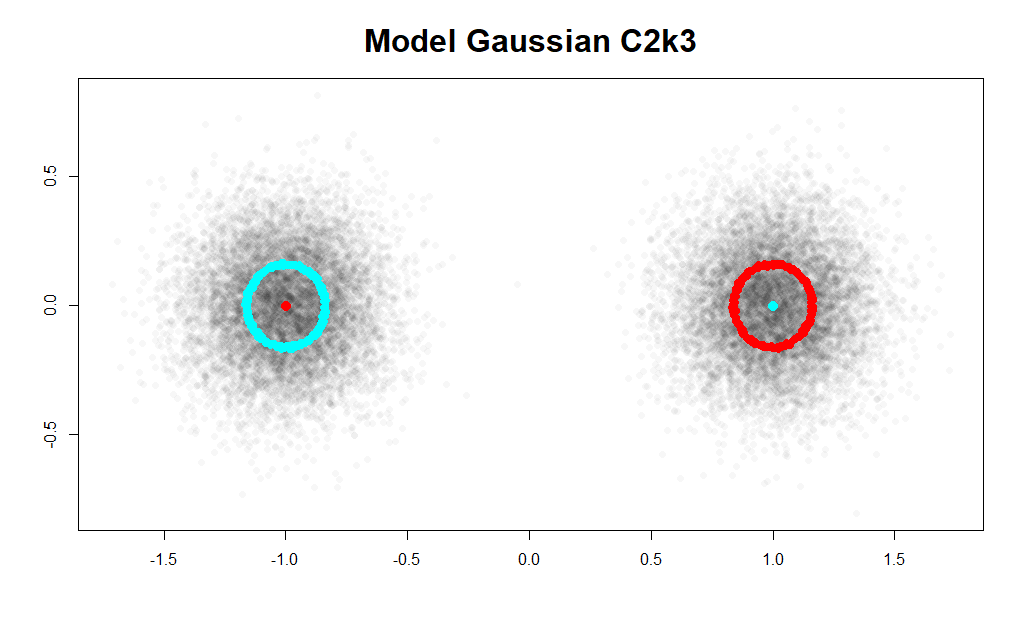} \\
    (a) & (b)
    \end{tabular}
    \caption{Empirical $k$-means of 500 simulated samples from model C1k2 (panel (a)) and C2k3 (panel (b)). The light gray points are $2\times 10^4$ observations from these distributions.
    }
    \label{fig:C1k2-C2k3-samples}
\end{figure}

\begin{figure}[!h]
				\centering
				\includegraphics[width=0.45\linewidth]{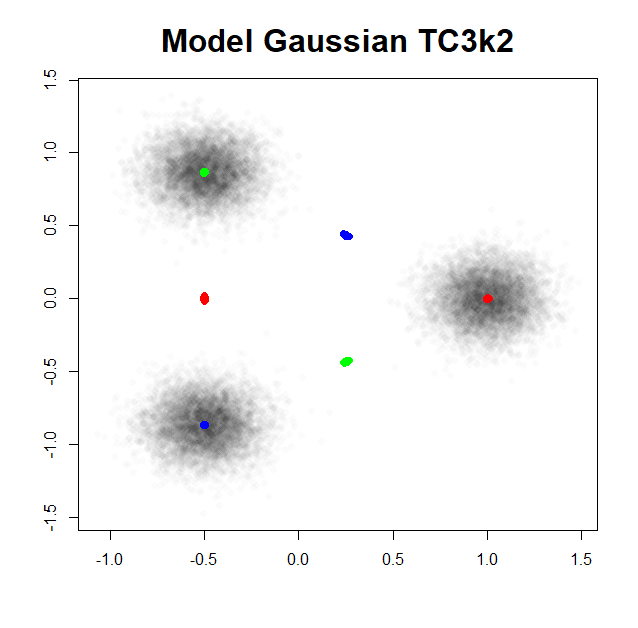}
				\caption{Empirical $k$-means of 500 simulated samples from model TC3k2. The light gray points are $2\times 10^4$ observations from this distribution. }
				\label{fig:TC3k2-samples}
\end{figure}

To illustrate the consistency result of $k$-means under non-uniqueness in Theorem \ref{th_cluster_consistency}, we have simulated 500 samples of size $n=2\times10^{4}$ from models C1k2, C2k3 and TC3k2 in Examples \ref{Example-C1k2}, \ref{ExampleC2k3} and \ref{Example-TC3k2}, respectively. For each of these samples, we have run the $k$-means algorithm and found the resulting set of centers. Therefore, we have calculated the empirical estimators of the $k$-means sets from large samples of the respective models. Figure \ref{fig:C1k2-C2k3-samples} (a) shows in red the empirical $k$-means sets of the 500 samples from model C1k2. We observe how the 2-centers of the samples completely cover the circumference in white which is the set of all population $k$-means. A completely analogous situation can be seen in Figure \ref{fig:C1k2-C2k3-samples} (b). Figure \ref{fig:TC3k2-samples} shows the results for the model TC3k2, which corresponds to discrete non-uniqueness.

\subsection{Practical conclusions}

 A few practical, preliminary, conclusions can be drawn from our empirical results.
\begin{enumerate}
    \item [(1)] The overall results are consistent with our theory in the sense that the asymptotic test for uniqueness seems to work, regarding both, the control of Type I error and the ability of \color{black}  detecting departures from the null hypothesis. Still, large sample sizes are needed in order to control the Type I error.

    \item [(2)] In the uniqueness example with a ``wrong'' choice of $k$ (model C3k2), we need an even larger sample size to control the significance level. This is not necessarily a bad feature since the result of the uniqueness test can be used to rule out values of $k$ that generate instability in the algorithm. Clearly, further study would be necessary to analyze how to use the uniqueness test to choose a reasonable value for the parameter $k$.

    \item [(3)] When we add noise variables to our population the problem becomes more complicated and the sample size required for satisfactory test performance is larger.

    \item [(4)] The analysis of the execution times of the $k$-means algorithm clearly shows that there is a distinction between the two cases of non-uniqueness that we have introduced in this work. In the case of continuous non-uniqueness the processing time is significantly longer than in the  discrete one. In the discrete non-uniqueness model TC3k2 the running-times are comparable to those of uniqueness cases with a ``correct'' selection of $k$. These empirical findings suggest  that the equivalence between stability and uniqueness repeatedly mentioned in the machine learning literature is not entirely accurate.

    \item [(5)] As an additional application, the ideas developed in this work  provide a numerical methodology to determine, by means of Monte Carlo simulations, whether a given (totally known) probability distribution satisfies the uniqueness property of the set of $k$-means, which might be a difficult task to perform by direct, ``analytical'' methods.

 \end{enumerate}


\section*{Acknowledgements}

We are very grateful to Professor Juan Cuesta-Albertos (Universidad de Cantabria) who read a preliminary version of this paper and made several enlightening comments. In particular, he pointed out to us the repulsion effect of $k$-means for distributions with support in the whole space. We also thank the help from Francisco de la Hoz (University of the Basque Country) with the computations.
This research has been partially supported by Grants PID2019-109387GB-I00, PID2023-148081NB-I00 from the Spanish Ministry of Science and Innovation and (for the second author) Grant CEX2019-000904-S funded by MCIN/AEI/ 10.13039/501100011033.  L.A. Rodr\'{i}guez acknowledges funding from the German Science Foundation DFG Research Unit 5381: ``Mathematical Statistics in the Information Age–Statistical Efficiency and Computational Tractability''.

  \end{document}